\documentclass[a4paper,10pt]{article}
\usepackage[english]{babel}
\usepackage[latin1]{inputenc}
\usepackage{amssymb}
\usepackage{amsmath}
\usepackage{amscd}
\usepackage{latexsym}
\usepackage{amsthm}
\usepackage{eucal}
\usepackage{mathrsfs}
\usepackage{mathabx}
\theoremstyle{plain}
\newtheorem{thm}{Theorem}[section]
\newtheorem{cor}[thm]{Corollary}
\newtheorem{lem}[thm]{Lemma}
\newtheorem{prop}[thm]{Proposition}

\newtheorem{oss}[thm]{Remark}
\newtheorem{esem}[thm]{Example}
\theoremstyle{definition}
\newtheorem{defn}{Definition}[section]

\title{Factoriality properties of moduli spaces of sheaves on abelian and K3 surfaces}
\author{Arvid Perego, Antonio Rapagnetta}

\begin{document}

\maketitle

\begin{abstract}
In this paper we complete the determination of the index of factoriality of moduli spaces of semistable sheaves on an abelian or projective K3 surface $S$. If $v=2w$ is a Mukai vector, $w$ is primitive, $w^{2}=2$ and $H$ is a generic polarization, let $M_{v}(S,H)$ be the moduli space of $H-$semistable sheaves on $S$ with Mukai vector $v$. First, we describe in terms of $v$ the pure weight-two Hodge structure and the Beauville form on the second integral cohomology of the symplectic resolutions of $M_{v}(S,H)$ (when $S$ is K3) and of the fiber $K_{v}(S,H)$ of the Albanese map of $M_{v}(S,H)$ (when $S$ is abelian). Then, if $S$ is K3 we show that $M_{v}(S,H)$ is either locally factorial or $2-$factorial, and we give an example of both cases. If $S$ is abelian, we show that $M_{v}(S,H)$ and $K_{v}(S,H)$ are $2-$factorial.
\end{abstract}

\section{Introduction and notations}

The study of the local factoriality of moduli spaces of sheaves on smooth, projective varieties is a problem which appears naturally when one wants to study their Picard group. The question of local factoriality has been adressed and solved in some cases. 

In \cite{DN} Dr\'ezet and Narasimhan show that the moduli space $U_{C}(r,d)$ of semi- stable vector bundles of rank $r$ and degree $d$ on a smooth, projective curve $C$ of genus $g\geq 2$, is locally factorial. The same result is proved even for the moduli space $U_{C}(r,L)$ of semistable vector bundles or rank $r$ and determinant $L\in Pic(C)$. 

If $G$ is a semisimple group, the moduli space $M_{G}$ of semistable principal $G-$bundles on a smooth, projective curve is locally factorial if and only if $G$ is special in the sense of Serre: this is due to Beauville, Laszlo and Sorger for $G$ classical or $G=G_{2}$ (see \cite{LaS}, \cite{BLS}, \cite{La}), and to Sorger (see \cite{S}) and Boysal, Kumar (see \cite{BK}) for exceptional $G$. 

In \cite{D1} Dr\'ezet shows that the moduli space $M(r,c_{1},c_{2})$ of semistable sheaves on $\mathbb{P}^{2}$ of rank $r\geq 2$ and Chern classes $c_{1}$ and $c_{2}$ is locally factorial. In \cite{D2}, Dr\'ezet studies the moduli space $M_{H}(r,c_{1},c_{2})$ of $H-$semistable sheaves on a rational, ruled surface $S$: he conjectures that it is locally factorial whenever the polarization $H$ is generic, otherwise he presents an explicit example of non locally factorial moduli space. In \cite{Y1} and \cite{Y2}, Yoshioka shows that the moduli space $M_{H}(r,c_{1},c_{2})$ of semistable sheaves on a ruled, non-rational surface is locally factorial whenever $H$ is generic. 

More recently, Kaledin, Lehn and Sorger in \cite{KLS} showed the local factoriality of a large class of moduli spaces of sheaves over projective K3 or abelian surfaces: in their work, the local factoriality plays a major role in order to avoid the existence of symplectic resolutions of the moduli space. 

The aim of this paper is to study the local factoriality of the remaing cases of moduli spaces of $H-$semistable sheaves on an abelian or projective K3 surface with generic polarization. In the following, $S$ will denote an abelian or projective K3 surface. An element $v\in \widetilde{H}(S,\mathbb{Z}):=H^{2*}(S,\mathbb{Z})$ will be written as $v=(v_{0},v_{1},v_{2})$, where $v_{i}\in H^{2i}(S,\mathbb{Z})$, and $v_{0},v_{2}\in\mathbb{Z}$. It will be called \textit{Mukai vector} if $v_{0}\geq 0$ and $v_{1}\in NS(S)$, and if $v_{0}=0$ then either $v_{1}$ is the first Chern class of an effective divisor, or $v_{1}=0$ and $v_{2}>0$. Recall that $\widetilde{H}(S,\mathbb{Z})$ has a pure weight-two Hodge structure defined as $$\widetilde{H}^{2,0}(S):=H^{2,0}(S),\,\,\,\,\,\,\,\,\,\,\,\widetilde{H}^{0,2}(S):=H^{0,2}(S),$$ $$\widetilde{H}^{1,1}(S):=H^{0}(S,\mathbb{C})\oplus H^{1,1}(S)\oplus H^{4}(S,\mathbb{C}),$$and a compatible lattice structure with respect to the Mukai pairing $(.,.)$. In the following, we let $v^{2}:=(v,v)$ for every Mukai vector $v$, and we define the sublattice $$v^{\perp}:=\{\alpha\in \widetilde{H}(S,\mathbb{Z})\,|\,(\alpha,v)=0\}\subseteq\widetilde{H}(S,\mathbb{Z}),$$which inherits a pure weight-two Hodge structure from the one on $\widetilde{H}(S,\mathbb{Z})$.

If $\mathscr{F}$ is a coherent sheaf on $S$, its \textit{Mukai vector} is $v(\mathscr{F}):=ch(\mathscr{F})\sqrt{td(S)}$. Let $H$ be a polarization and $v$ a Mukai vector on $S$. We write $M_{v}(S,H)$ (resp. $M_{v}^{s}(S,H)$) for the moduli space of $H-$semistable (resp. $H-$stable) sheaves on $S$ with Mukai vector $v$. If no confusion on $S$ and $H$ is possible, we drop them from the notation.

From now on, we suppose that $H$ is a $v-$generic polarization (for a definition, see \cite{PR}). We write $v=mw$, where $m\in\mathbb{N}$ and $w$ is a primitive Mukai vector on $S$. It is known that if $M_{v}^{s}\neq\emptyset$, then $M_{v}^{s}$ is smooth, quasi-projective, of dimension $v^{2}+2$ and carries a symplectic form (see Mukai \cite{M1}). 

If $S$ is an abelian surface, a further construction is necessary. Choose $\mathscr{F}_{0}\in M_{v}(S,H)$, and define $a_{v}:M_{v}(S,H)\longrightarrow S\times\widehat{S}$ in the following way (see \cite{Y4}): let $p_{\widehat{S}}:S\times \widehat{S}\longrightarrow\widehat{S}$ be the projection and $\mathscr{P}$ the Poincar\'e bundle on $S\times\widehat{S}$. For every $\mathscr{F}\in M_{v}(S,H)$ we let $$a_{v}(\mathscr{F}):=(det(p_{\widehat{S}!}((\mathscr{F}-\mathscr{F}_{0})\otimes(\mathscr{P}-\mathscr{O}_{S\times\widehat{S}})),det(\mathscr{F})\otimes det(\mathscr{F}_{0})^{-1}).$$Moreover, we define $K_{v}(S,H):=a_{v}^{-1}(0_{S},\mathscr{O}_{S})$, where $0_{S}$ is the zero of $S$.

The local factoriality of $M_{v}(S,H)$ and $K_{v}(S,H)$ is almost completely understood. First, write $v=mw$, and suppose that $w^{2}\geq 0$, as if $w^{2}<0$ then $M_{v}(S,H)$ is either empty or a single point. Under this condition, then $M_{v}(S,H)$ and $K_{v}(S,H)$ are known to be locally factorial in one of the two following cases (see the Appendix for the details on $K_{v}$):
\begin{enumerate}
 \item $m=1$ and $w^{2}\geq 0$ (in this case $M_{v}=M^{s}_{v}$, hence it is smooth);
 \item $m=2$ and $w^{2}\geq 4$, or $m\geq 3$ and $w^{2}\geq 2$ (see \cite{KLS}).
\end{enumerate}

The aim of the present paper is to determine the factoriality index of the remaining cases, namely when $m\geq 2$ and $w^{2}=0$, and when $m=2$ and $w^{2}=2$. If $m\geq 2$ and $w^{2}=0$, then $M_{w}$ is either a projective K3 surface or an abelian surface, and $M_{v}$ is isomorphic to the $m-$th symmetric product $M_{w}^{(m)}$. Hence, it is $2-$factorial: as we did not find any reference for this, we included the proof in section 2. 

Sections 3 and 4 are devoted to the factorialiy index of the last case: $m=2$ and $w^{2}=2$, i. e. $(S,v,H)$ is an OLS-triple following the terminology of \cite{PR}. A particular case is given by the moduli spaces $M_{10}$ and $K_{6}$ described by O'Grady in \cite{OG1} and \cite{OG2}, which are known to be $2-$factorial by \cite{Pe}. We recall that if $(S,v,H)$ is on OLS-triple, then $M_{v}$ (resp. $K_{v}$) admits a symplectic resolution $\widetilde{M}_{v}$ (resp. $\widetilde{K}_{v}$) (see \cite{LS}). In \cite{PR} we proved that $\widetilde{M}_{v}$ and $\widetilde{K}_{v}$ are irreducible symplectic manifolds, and we showed that $H^{2}(M_{v},\mathbb{Z})$ and $H^{2}(K_{v},\mathbb{Z})$ have a pure weight-two Hodge structure, a lattice structure, and they are Hodge isometric to $v^{\perp}$. 

As a first step to compute the index of factoriality of $M_{v}$ and $K_{v}$, we use the results of \cite{PR} to calculate the pure weight-two Hodge structure and the lattice structure on $H^{2}(\widetilde{M}_{v},\mathbb{Z})$ (if $S$ is K3) and on $H^{2}(\widetilde{K}_{v},\mathbb{Z})$ in terms of the Mukai vector $v$. We show that if $S$ is K3, then $H^{2}(\widetilde{M}_{v},\mathbb{Z})$ is Hodge isometric to a $\mathbb{Z}-$submodule $\Gamma_{v}$ of $(v^{\perp}\otimes\mathbb{Q})\oplus\mathbb{Q}$, on which we define a pure weight-two Hodge structure and a lattice structure coming from those on $v^{\perp}$. This is the content of Theorem \ref{thm:w2hs}, where a similar description of $H^{2}(\widetilde{K}_{v},\mathbb{Z})$ is presented. These results are analogues to the description of the Hodge structure and of the lattice structure on $H^{2}(M_{v},\mathbb{Z})$ and $H^{2}(K_{v},\mathbb{Z})$ given by O'Grady and Yoshioka for primitive $v$.

Once the description of the Hodge structure on $H^{2}(\widetilde{M}_{v},\mathbb{Z})$ and $H^{2}(\widetilde{K}_{v},\mathbb{Z})$ is done, we proceed with the determination of the index of factoriality of $M_{v}(S,H)$ and $K_{v}(S,H)$. More precisely, if $S$ is a K3 surface we show the following:

\begin{thm}
\label{thm:maink3}
Let $(S,v,H)$ be an OLS-triple such that $S$ is a projective K3 surface, and write $v=2w$. Then $M_{v}(S,H)$ is either $2-$factorial or  locally factorial. Moreover
 \begin{enumerate}
  \item $M_{v}(S,H)$ is $2-$factorial if and only if there is $\gamma\in \widetilde{H}(S,\mathbb{Z})\cap\widetilde{H}^{1,1}(S)$ such that $(\gamma,w)=1$.
  \item $M_{v}(S,H)$ is locally factorial if and only if for every $\gamma\in\widetilde{H}(S,\mathbb{Z})\cap\widetilde{H}^{1,1}(S)$ we have $(\gamma,w)\in 2\mathbb{Z}$.
 \end{enumerate}
\end{thm}

To prove this result, we use the description of the Hodge structure on $H^{2}(M_{v},\mathbb{Z})$ given in \cite{PR}, and of the Hodge structure of $H^{2}(\widetilde{M}_{v},\mathbb{Z})$ we obtain in Theorem \ref{thm:w2hs}: we are then able to calculate the Picard groups of $\widetilde{M}_{v}$ and $M_{v}$, and from the relation between them we determine the factoriality index. This Theorem allows us to easily present explicit examples of OLS-triples $(S,v,H)$ such that $M_{v}(S,H)$ is locally factorial (see Example \ref{esem:lfmv}).

In the case of OLS-triples $(S,v,H)$ where $S$ is abelian, we study the local factoriality of both $M_{v}(S,H)$ and $K_{v}(S,H)$. The result we show is the following:

\begin{thm}
\label{thm:mainab}
Let $(S,v,H)$ be an OLS-triple such that $S$ is an abelian surface. Then $M_{v}(S,H)$ and $K_{v}(S,H)$ are $2-$factorial.
\end{thm}

We notice that this result does not depend on $v$: the $2-$factoriality of $K_{v}$ is a consequence of the existence of a square root of the exceptional divisor $\widetilde{\Sigma}_{v}$ of the symplectic resolution $\widetilde{K}_{v}\longrightarrow K_{v}$. The $2-$factoriality of $M_{v}$ is more involved and requires a different argument (see section 4.2.2), as here the exceptional divisor of the symplectic resolution $\widetilde{M}_{v}\longrightarrow M_{v}$ does not admit any square root.

\section{Recall on factoriality and on the O'Grady examples}

In this section we recall some basic elements about the notions of local factoriality and of $2-$factoriality, and we resume some basic results about the O'Grady examples introduced in \cite{OG1} and \cite{OG2} that will be used in the following. Moreover, we recall the definition of OLS-triple and the main results of \cite{PR}.

\subsection{Factoriality}

Let $X$ be a normal, projective variety.

\begin{defn}
We say that $X$ is \textit{locally factorial} if for every $x\in X$ the ring $\mathscr{O}_{X,x}$ is a unique factorization domain.
\end{defn}

We denote by $A^{1}(X)$ the group of linear equivalence classes of Weil divisors on $X$. If $D$ is a Weil divisor on $X$, we still write $D$ for its linear equivalence class. We have a natural inclusion $d_{X}:Pic(X)\longrightarrow A^{1}(X)$ associating to any line bundle (or Cartier divisor) the corresponding Weil divisor.

\begin{defn}
We say that $X$ is \textit{factorial} if the morpism $d_{X}$ is an isomorphism, and $n-$\textit{factorial} if co$\ker(d_{X})$ is $n-$torsion.
\end{defn}

It is well-known that $X$ is locally factorial if and only if it is factorial. 

\begin{oss}
\label{oss:cw}{\rm We recall a basic fact which will be used several times in the following. Let $U\subseteq X$ be an open subset, and let $Z:=X\setminus U$. If $codim_{X}(Z)\geq 2$, then we have $A^{1}(X)=A^{1}(U)$. If $U$ is smooth, then we have an identification $A^{1}(X)=Pic(U)$. A particular case is when $U$ is the smooth locus of $X$, as $X$ is normal. In this case, we have $Pic(X)\subseteq Pic(U)$: it follows that if $\pi:\widetilde{X}\longrightarrow X$ is a resolution of the singularities, then $\pi^{*}:Pic(X)\longrightarrow Pic(\widetilde{X})$ is injective.}
\end{oss}

As we recalled in the introduction, if $S$ is an abelian or K3 surface, $v=mw$ is a Mukai vector on $S$ such that $w$ is primitive and $w^{2}\geq 0$, and $H$ is a $v-$generic polarization, the moduli spaces $M_{v}(S,H)$ and $K_{v}(S,H)$ are known to be locally factorial if $m=1$, if $m=2$ and $w^{2}\geq 4$, and if $m\geq 3$ and $w^{2}\geq 2$. The aim of this paper is to study the index of factoriality of the remaining cases, namely $m\geq 2$ and $w^{2}=0$, or $m=2$ and $w^{2}=2$.

If $m\geq 2$ and $w^{2}=0$, the moduli space $M_{v}(S,H)$ is isomorphic to $M_{w}^{(m)}$, the $m-$th symmetric product of $M_{w}$, which is either a K3 surface (if $S$ is K3) or an abelian surface (if $S$ is abelian). It follows that $M_{v}$ is $2-$factorial: as we did not find any reference for this, we included the proof in the following example.

\begin{esem}
\label{esem:symmk3}
{\rm Let $S$ be a smooth projective complex surface, and $n\in\mathbb{N}$, $n\geq 2$. The $n-$th symmetric product $S^{(n)}$ is singular, and its singular locus $\Delta$ has codimension 2 in $S^{(n)}$. Let $Hilb^{n}(S)$ be the Hilbert scheme of $n-$points on $S$, and let $\rho_{n}:Hilb^{n}(S)\longrightarrow S^{(n)}$ be the Hilbert-Chow morphism. By \cite{F}, we have that $\rho_{n}$ is a resolution of the singularities, and we write $E$ for the exceptional divisor, which is irreducible.

To show that $S^{(n)}$ is $2-$factorial, we consider the following commutative diagram:
$$\begin{array}{ccc}
Pic(S^{(n)}) & \stackrel{\rho_{n}^{*}}\longrightarrow & Pic(Hilb^{n}(S))\\
\scriptstyle{j}\downarrow & & \downarrow\scriptstyle{g}\\
Pic(S^{n})^{\sigma_{n}} & \stackrel{f}\longrightarrow & Pic(S^{n})^{\sigma_{n}}\oplus\mathbb{Z}
\end{array}$$where $Pic(S^{n})^{\sigma_{n}}$ is the group of the line bundles on $S^{n}$ which are invariant under the action of the symmetric group $\sigma_{n}$, the morphism $f$ is the identity on $Pic(S^{n})^{\sigma_{n}}$, and $j$ and $g$ are two natural morphisms (for a definition, see \cite{F}). The morphism $\rho_{n}^{*}$ is injective (see Remark \ref{oss:cw}), the morphism $j$ is surjective by Lemma 6.1 of \cite{F}, and by Theorem 6.2 of \cite{F} $g$ is an isomorphism such that $g(\mathscr{O}(E))=(0,2)$. From this, it follows that 
\begin{equation}
\label{eq:pichilb}
Pic(Hilb^{n}(S))=\rho_{n}^{*}(Pic(S^{(n)}))\oplus\mathbb{Z}\cdot A,
\end{equation}
where $A\in Pic(Hilb^{n}(S))$ is a line bundle such that $A^{\otimes 2}=\mathscr{O}(E)$. Now, let $i:Hilb^{n}(S)\setminus E\longrightarrow Hilb^{n}(S)$ be the inclusion, and consider the following exact sequence (see Proposition 6.5 in Chapter II of \cite{H}) $$0\longrightarrow\mathbb{Z}\stackrel{\iota}\longrightarrow Pic(Hilb^{n}(S))\stackrel{i^{*}}\longrightarrow Pic(Hilb^{n}(S)\setminus E)\longrightarrow 0,$$where $\iota(1):=\mathscr{O}(E)$, so that $$Pic(Hilb^{n}(S)\setminus E)=Pic(Hilb^{n}(S))/\mathbb{Z}\cdot\mathscr{O}(E)=\rho_{n}^{*}(Pic(S^{(n)}))\oplus\mathbb{Z}/2\mathbb{Z}\cdot A.$$By Remark \ref{oss:cw} we have $A^{1}(S^{(n)})=Pic(Hilb^{n}(S)\setminus E)$, so that equation (\ref{eq:pichilb}) gives $$A^{1}(S^{(n)})/Pic(S^{(n)})=(\rho_{n}^{*}(Pic(S^{(n)}))\oplus\mathbb{Z}/2\mathbb{Z}\cdot A)/\rho_{n}^{*}(Pic(S^{(n)}))\simeq \mathbb{Z}/2\mathbb{Z},$$and the $2-$factoriality of $S^{(n)}$ is shown.}
\end{esem}

\subsection{Recall on the O'Grady examples}

We now recall some basic elements of the two O'Grady examples. We start with the $10-$dimensional one: let $X$ be a projective K3 surface such that $Pic(X)=\mathbb{Z}\cdot H$, where $H$ is ample and $H^{2}=2$, and let $v:=(2,0,-2)$ be a Mukai vector on $X$. We let $M_{10}:=M_{v}(X,H)$, which is $10-$dimensional, its smooth locus is $M_{10}^{s}:=M_{v}^{s}(X,H)$, and its singular locus is denoted $\Sigma$. In \cite{OG1} O'Grady shows that $M_{10}$ admits a symplectic resolution $\pi:\widetilde{M}_{10}\longrightarrow M_{10}$ and that $\widetilde{M}_{10}$ is an irreducible symplectic manifold. Here are the main results we will use in the following about $\widetilde{M}_{10}$ (see \cite{R1} and \cite{Pe}):

\begin{thm}
\label{thm:m10}Let $B\subseteq M_{10}$ be the Weil divisor which parameterizes non-locally free sheaves, $\widetilde{\Sigma}$ the exceptional divisor of $\pi$ and $\widetilde{B}$ the proper transform of $B$ under $\pi$.
\begin{enumerate}
 \item We have $$H^{2}(\widetilde{M}_{10},\mathbb{Z})=\widetilde{\mu}(H^{2}(X,\mathbb{Z}))\oplus\mathbb{Z}\cdot c_{1}(\widetilde{\Sigma})\oplus\mathbb{Z}\cdot c_{1}(\widetilde{B}),$$where $\widetilde{\mu}:H^{2}(X,\mathbb{Z})\longrightarrow H^{2}(\widetilde{M}_{10},\mathbb{Z})$ is the Donaldson morphism.
 \item Let $B_{10}$ be the Beauville form of $\widetilde{M}_{10}$. For every $\alpha\in H^{2}(X,\mathbb{Z})$ we have $$B_{10}(\widetilde{\mu}(\alpha),\widetilde{\Sigma})=0,\,\,\,\,\,\,\,\,\,\,\,B_{10}(\widetilde{\mu}(\alpha),\widetilde{B})=0,$$ $$B_{10}(\widetilde{\Sigma},\widetilde{\Sigma})=-6,\,\,\,\,\,\,\,\,\,\,\,B_{10}(\widetilde{\Sigma},\widetilde{B})=3,\,\,\,\,\,\,\,\,\,\,B_{10}(\widetilde{B},\widetilde{B})=-2.$$
 \item We have $v^{\perp}=\{(n,\alpha,n)\,|\,n\in\mathbb{Z},\,\,\alpha\in H^{2}(X,\mathbb{Z})\}$, and the morphism $$\widetilde{\lambda}:v^{\perp}\longrightarrow H^{2}(\widetilde{M}_{10},\mathbb{Z}),\,\,\,\,\,\,\,\,\widetilde{\lambda}(n,\alpha,n)=\widetilde{\mu}(\alpha)+2nc_{1}(\widetilde{B})+nc_{1}(\widetilde{\Sigma})$$is a Hodge isometry onto its image.
 \item The moduli space $M_{10}$ is $2-$factorial.
 \end{enumerate}
\end{thm}

For the $6-$dimensional O'Grady example, let $X$ be an abelian surface such that $NS(X)=\mathbb{Z}\cdot H$, where $H$ is ample and $H^{2}=2$, and let $v:=(2,0,-2)$ be a Mukai vector on $X$. We let $K_{6}:=K_{v}(X,H)$, which is $6-$dimensional, its smooth locus is $K_{6}^{s}:=K_{v}^{s}(X,H)$, and its singular locus is denoted $\Sigma$. In \cite{OG2} O'Grady shows that $K_{6}$ admits a symplectic resolution $\pi:\widetilde{K}_{6}\longrightarrow K_{6}$ and that $\widetilde{K}_{6}$ is irreducible symplectic. We have the following (see \cite{R2} and \cite{Pe}):

\begin{thm}
\label{thm:k6}Let $B\subseteq K_{6}$ be the Weil divisor which parameterizes non-locally free sheaves, $\widetilde{\Sigma}$ the exceptional divisor of $\pi$ and $\widetilde{B}$ the proper transform of $B$ under $\pi$.
\begin{enumerate}
 \item We have $$H^{2}(\widetilde{K}_{6},\mathbb{Z})=\widetilde{\mu}(H^{2}(X,\mathbb{Z}))\oplus\mathbb{Z}\cdot c_{1}(A)\oplus\mathbb{Z}\cdot c_{1}(\widetilde{B}),$$where $\widetilde{\mu}:H^{2}(X,\mathbb{Z})\longrightarrow H^{2}(\widetilde{K}_{6},\mathbb{Z})$ is the Donaldson morphism, and $A\in Pic(\widetilde{K}_{6})$ is such that $A^{\otimes 2}=\mathscr{O}(\widetilde{\Sigma})$.
 \item Let $B_{6}$ be the Beauville form of $\widetilde{K}_{6}$. Then for every $\alpha\in H^{2}(X,\mathbb{Z})$ we have $$B_{6}(\widetilde{\mu}(\alpha),A)=0,\,\,\,\,\,\,\,\,\,\,\,B_{6}(\widetilde{\mu}(\alpha),\widetilde{B})=0,$$ $$B_{6}(A,A)=-2,\,\,\,\,\,\,\,\,\,\,\,B_{6}(A,\widetilde{B})=2,\,\,\,\,\,\,\,\,\,\,B_{6}(\widetilde{B},\widetilde{B})=-4.$$
 \item We have $v^{\perp}=\{(n,\alpha,n)\,|\,n\in\mathbb{Z},\,\,\alpha\in H^{2}(X,\mathbb{Z})\}$, and the morphism $$\widetilde{\lambda}:v^{\perp}\longrightarrow H^{2}(\widetilde{K}_{6},\mathbb{Z}),\,\,\,\,\,\,\,\,\widetilde{\lambda}(n,\alpha,n)=\widetilde{\mu}(\alpha)+nc_{1}(\widetilde{B})+nc_{1}(A)$$is a Hodge isometry onto its image.
 \item The moduli space $K_{6}$ is $2-$factorial.
 \end{enumerate}
\end{thm}

The two O'Grady examples are generalized with the following definition, which is contained in \cite{PR}:

\begin{defn}
Let $S$ be an abelian or projective K3 surface, $v$ a Mukai vector, $H$ an ample line bundle on $S$. We say that $(S,v,H)$ is an \textit{OLS-triple} if the following conditions are verified:
\begin{enumerate}
 \item the polarization $H$ is primitive and $v-$generic;
 \item there is a primitive Mukai vector $w\in\widetilde{H}(S,\mathbb{Z})$ such that $v=2w$ and $w^{2}=2$;
 \item if $w=(0,\xi,a)$, then $a\neq 0$.
\end{enumerate}
\end{defn}

If $(S,v,H)$ is an OLS-triple, then $M_{v}(S,H)$ is a normal, projective variety of dimension 10, whose smooth locus is $M_{v}^{s}(S,H)$, and whose singular locus is denoted $\Sigma_{v}$, which has codimension 2 in $M_{v}$. By \cite{LS}, we know that $M_{v}$ admits a symplectic resolution $\pi_{v}:\widetilde{M}_{v}(S,H)\longrightarrow M_{v}(S,H)$ obtained by blowing up $M_{v}$ along $\Sigma_{v}$ with reduced structure. If $S$ is abelian, let $\widetilde{K}_{v}(S,H):=\pi_{v}^{-1}(K_{v}(S,H))$, which is a symplectic resolution of $K_{v}(S,H)$. By abuse of notation, we still write $\pi_{v}:\widetilde{K}_{v}\longrightarrow K_{v}$. We have the following result, which is Theorems 1.3 and 1.5 of \cite{PR}.

\begin{thm}
\label{thm:pr}Let $(S,v,H)$ be an OLS-triple.
\begin{enumerate}
 \item If $S$ is K3, then $\widetilde{M}_{v}$ is irreducible symplectic and deformation equivalent to $\widetilde{M}_{10}$. The morphism $\pi_{v}^{*}:H^{2}(M_{v},\mathbb{Z})\longrightarrow H^{2}(\widetilde{M}_{v},\mathbb{Z})$ is injective, and the restrictions to $H^{2}(M_{v},\mathbb{Z})$ of the pure weight-two Hodge structure and of the Beauville form on $H^{2}(\widetilde{M}_{v},\mathbb{Z})$ give a pure weight-two Hodge structure and a compatible lattice structure on $H^{2}(M_{v},\mathbb{Z})$. Moreover, there is a Hodge isometry $\lambda_{v}:v^{\perp}\longrightarrow H^{2}(M_{v},\mathbb{Z})$.
 \item If $S$ is abelian, then $\widetilde{K}_{v}$ is irreducible symplectic and deformation equivalent to $\widetilde{K}_{6}$. The morphism $\pi_{v}^{*}:H^{2}(K_{v},\mathbb{Z})\longrightarrow H^{2}(\widetilde{K}_{v},\mathbb{Z})$ is injective, and the restrictions to $H^{2}(K_{v},\mathbb{Z})$ of the pure weight-two Hodge structure and of the Beauville form on $H^{2}(\widetilde{K}_{v},\mathbb{Z})$ give a pure weight-two Hodge structure and a compatible lattice structure on $H^{2}(K_{v},\mathbb{Z})$. Moreover, there is a Hodge isometry $\nu_{v}:v^{\perp}\longrightarrow H^{2}(K_{v},\mathbb{Z})$.
\end{enumerate}
\end{thm}

The morphisms $\lambda_{v}$ and $\nu_{v}$ are constructed in section 3.2 of \cite{PR}.

\begin{oss}
\label{oss:lambdatilde}{\rm If $(S,v,H)$ is an OLS-triple and $S$ is K3, it follows from point 1 of Theorem \ref{thm:pr} that the morphism $$\widetilde{\lambda}_{v}:=\pi_{v}^{*}\circ\lambda_{v}:v^{\perp}\longrightarrow H^{2}(\widetilde{M}_{v},\mathbb{Z})$$is an injective Hodge morphism which is an isometry onto its image. If $S$ is abelian, we still have a morphism $\lambda_{v}:v^{\perp}\longrightarrow H^{2}(M_{v},\mathbb{Z})$, which is such that $\nu_{v}=j_{v}^{*}\circ\lambda_{v}$, where $j_{v}:K_{v}\longrightarrow M_{v}$ is the inclusion. We define $\widetilde{\lambda}_{v}:=\pi_{v}^{*}\circ\lambda_{v}$ and $$\widetilde{\nu}_{v}:=\pi_{v}^{*}\circ\nu_{v}:v^{\perp}\longrightarrow H^{2}(\widetilde{K}_{v},\mathbb{Z}).$$As a consequence of point 2 of Theorem \ref{thm:pr}, we have that $\widetilde{\nu}_{v}$ is an injective Hodge morphism which is an isometry ontro its image. We even remark that $\widetilde{\nu}_{v}=\widetilde{j}_{v}^{*}\circ\widetilde{\lambda}_{v}$, where $\widetilde{j}_{v}:\widetilde{K}_{v}\longrightarrow\widetilde{M}_{v}$ is the inclusion.}
\end{oss}

As a consequence of Theorem \ref{thm:pr} we calculate the Picard groups of $\widetilde{M}_{v}$ (if $S$ is K3) and of $\widetilde{K}_{v}$ (if $S$ is abelian):

\begin{cor}
\label{cor:picmv}Let $(S,v,H)$ be an OLS-triple.
\begin{enumerate}
 \item If $S$ is K3, there is an isomorphism $L_{v}:v^{\perp}\cap\widetilde{H}^{1,1}(S)\longrightarrow Pic(M_{v})$ such that for every $\alpha\in v^{\perp}\cap\widetilde{H}^{1,1}(S)$ we have $\lambda_{v}(\alpha)=c_{1}(L_{v}(\alpha))$.
 \item If $S$ is abelian, there is an isomorphism $N_{v}:v^{\perp}\cap\widetilde{H}^{1,1}(S)\longrightarrow Pic(K_{v})$, such that for every $\alpha\in v^{\perp}\cap\widetilde{H}^{1,1}(S)$ we have $\nu_{v}(\alpha)=c_{1}(N_{v}(\alpha))$.
\end{enumerate}
\end{cor}

\proof If $S$ is K3, the morphism $L_{v}$ is constructed in section 3.2 of \cite{PR} by using Le Potier's determinant construction, and for every $\alpha\in v^{\perp}\cap\widetilde{H}^{1,1}(S)$ we have $\lambda_{v}(\alpha)=c_{1}(L_{v}(\alpha))$ by definition. By point 1 of Theorem \ref{thm:pr} we know that $\lambda_{v}$ is an isomorphism: this implies immediately that $L_{v}$ is injective. For the surjectivity, if $L\in Pic(M_{v})$, then $\pi_{v}^{*}(c_{1}(L))\in im(\widetilde{\lambda}_{v})\cap H^{1,1}(\widetilde{M}_{v})$, hence there is $\alpha\in v^{\perp}\cap\widetilde{H}^{1,1}(S)$ such that $\pi_{v}^{*}(c_{1}(L))=\widetilde{\lambda}_{v}(\alpha)$. As $\widetilde{M}_{v}$ is irreducible symplectic, it follows that $\pi_{v}^{*}(L)$ and $\pi_{v}^{*}(L_{v}(\alpha))$ are two isomorphic line bundles. Now, $\pi_{v}^{*}$ injects $Pic(M_{v})$ in $Pic(\widetilde{M}_{v})$ (see Remark \ref{oss:cw}), hence $L$ and $L_{v}(\alpha)$ are isomorphic, and the surjectivity of $L_{v}$ is shown.

If $S$ is abelian, Le Potier's determinant construction gives us a morphism $L_{v}:e_{v}^{\perp}\longrightarrow Pic(K_{v})$, where $e_{v}\in K_{hol}(S)$ is the class of a sheaf parameterized by $M_{v}$ in the homological $K-$theory of $S$, and the orthogonal is with respect to the natural pairing on $K_{hol}(S)$ given by the Euler characteristic. For every $\alpha\in v^{\perp}\cap \widetilde{H}^{1,1}(S)$, we then define $N_{v}(\alpha):=j_{v}^{*}(L_{v}(F))$, where $F\in K_{hol}(S)$ is an element whose Chern character is $\alpha$, and $j_{v}:K_{v}\longrightarrow M_{v}$ is the inclusion. The proof then goes as in the previous case, where one uses $K_{v}$ instead of $M_{v}$, $N_{v}$ instead of $L_{v}$, $j_{v}^{*}\circ\lambda_{v}$ instead of $\lambda_{v}$, and point 2 of Theorem \ref{thm:pr} instead of point 1.\endproof

\section{The weight-two Hodge structure of the symplectic resolutions}

In this section we calculate the weight-two Hodge structures on $H^{2}(\widetilde{M}_{v},\mathbb{Z})$ (when $S$ is K3) and on $H^{2}(\widetilde{K}_{v},\mathbb{Z})$ (when $S$ is abelian), together with their Beauville form. This result is interesting even on its own, and it will be used in the remaining part of the paper to prove Theorems \ref{thm:maink3} and \ref{thm:mainab}.

The first step to describe the weight-two Hodge structures on $H^{2}(\widetilde{M}_{v},\mathbb{Z})$ and on $H^{2}(\widetilde{K}_{v},\mathbb{Z})$ is to describe the lattice structures (with respect to the Beauville form) in terms of the Mukai vector $v$. In the following, we will make use of the notation $\oplus_{\perp}$ for the orthogonal sum of lattices, using it even for orthogonal direct sums of lattices whose bilinear symmetric form takes values in $\mathbb{Q}$. The description of $H^{2}(\widetilde{K}_{v},\mathbb{Z})$ is easy (we will se that in this case $H^{2}(\widetilde{K}_{v},\mathbb{Z})$ is Hodge isometric to $v^{\perp}\oplus_{\perp}\mathbb{Z}(-2)$, where $\mathbb{Z}(-2)$ is the rank 1 lattice generated by an element of square $-2$), but more involved in the case of K3 surface. 

For this purpose, we need to introduce some notations: let $S$ be a projective K3 surface, and $w\in\widetilde{H}(S,\mathbb{Z})$ a Mukai vector such that $w^{2}=2$. Let $v:=2w$ and $\Delta(v^{\perp})$ the discriminant of the lattice $v^{\perp}$ (which is equal to $w^{\perp}$): we have $|\Delta(v^{\perp})|=w^{2}=2$. Let $(v^{\perp})^{*}$ be the dual lattice of $v^{\perp}$, i. e. $$(v^{\perp})^{*}=Hom_{\mathbb{Z}}(v^{\perp},\mathbb{Z})=\{\alpha\in v^{\perp}\otimes_{\mathbb{Z}}\mathbb{Q}\,|\,(\alpha,\beta)\in\mathbb{Z},\,\,\mathrm{for}\,\,\mathrm{every}\,\,\beta\in v^{\perp}\}.$$The bilinear symmetric form on $(v^{\perp})^{*}$ is the $\mathbb{Q}-$bilinear extension of the Mukai pairing on $v^{\perp}$, and we denote it by $(.,.)_{\mathbb{Q}}$: in general, it does not take values in $\mathbb{Z}$, but surely in $\mathbb{Q}$. 

The morphism $$v^{\perp}\longrightarrow(v^{\perp})^{*},\,\,\,\,\,\,\,\,\,\,\,\,\,\alpha\mapsto(\alpha,.):v^{\perp}\longrightarrow\mathbb{Z}$$is injective, and $v^{\perp}$ is a sublattice of $(v^{\perp})^{*}$ of index $|\Delta(v^{\perp})|=2$. Let us now consider the $\mathbb{Z}-$module $(v^{\perp})^{*}\oplus_{\perp}\mathbb{Z}\cdot(\sigma/2)$, over which we define a bilinear symmetric form $b_{v}$ by letting $b_{v|(v^{\perp})^{*}}:=(.,.)_{\mathbb{Q}}$ and $b_{v}(\sigma,\sigma):=-6$. This is a non-degenerate symmetric bilinear form taking values in $\mathbb{Q}$.

Now, let $\Gamma_{v}$ be the subset of $(v^{\perp})^{*}\oplus_{\perp}\mathbb{Z}\cdot(\sigma/2)$ defined as $$\Gamma_{v}:=\{(\alpha,k\sigma/2)\in(v^{\perp})^{*}\oplus_{\perp}\mathbb{Z}\cdot(\sigma/2)\,|\,k\in 2\mathbb{Z}\,\,\mathrm{if}\,\,\mathrm{and}\,\,\mathrm{only}\,\,\mathrm{if}\,\,\alpha\in v^{\perp}\}.$$It is easy to see that $\Gamma_{v}$ is a submodule of $(v^{\perp})^{*}\oplus_{\perp}\mathbb{Z}\cdot(\sigma/2)$.

\begin{oss}
\label{oss:explicitgamma}{\rm The submodule $\Gamma_{v}$ has a more explicit description as a submodule of $(v^{\perp}\otimes\mathbb{Q})\oplus_{\perp}\mathbb{Q}\cdot\sigma$. More precisely, we have $$\Gamma_{v}=\{(\beta/2,k\sigma/2)\,|\,\beta\in v^{\perp},\,\,k\in\mathbb{Z},\,\,(\beta,v^{\perp})\subseteq 2\mathbb{Z},\mathrm{and}\,\,k\in 2\mathbb{Z}\,\,\mathrm{iff}\,\,\beta/2\in v^{\perp}\}.$$Indeed, if $\beta\in v^{\perp}$ is such that $(\beta,v^{\perp})\subseteq 2\mathbb{Z}$, then $(\beta/2,v^{\perp})\subseteq\mathbb{Z}$, so that $\beta/2\in(v^{\perp})^{*}$. Conversely, if $\alpha\in(v^{\perp})^{*}$, then $\beta:=2\alpha\in v^{\perp}$ as the index of $v^{\perp}$ in $(v^{\perp})^{*}$ is 2. Hence $(\beta,v^{\perp})\subseteq 2\mathbb{Z}$. This proves that $$(v^{\perp})^{*}=\{\beta/2\in v^{\perp}\otimes\mathbb{Q}\,|\,\beta\in v^{\perp},\,\,\mathrm{and}\,\,(\beta,v^{\perp})\subseteq 2\mathbb{Z}\},$$and the equality above follows immediately.}
\end{oss}

\begin{oss}
\label{oss:isom10}{\rm If we choose $v=(2,0,-2)$, then it is easy to see that $\Gamma_{v}$ is isometric to $H^{2}(\widetilde{M}_{10},\mathbb{Z})$. Indeed, we have a $\mathbb{Q}-$linear morphism $$f_{v}:(v^{\perp}\otimes\mathbb{Q})\oplus_{\perp}\mathbb{Q}\cdot\sigma\longrightarrow H^{2}(\widetilde{M}_{10},\mathbb{Q}),\,\,\,\,\,\,\,\,\,\,f_{v}(\alpha,h\sigma):=\widetilde{\lambda}(\alpha)+hc_{1}(\widetilde{\Sigma}),$$where by abuse of notation we still write $\widetilde{\lambda}$ for the $\mathbb{Q}-$linear extension of $\widetilde{\lambda}$. Now, by definition of $b_{v}$ we have $$b_{v}((\alpha_{1},h_{1}\sigma),(\alpha_{2},h_{2}\sigma))=(\alpha_{1},\alpha_{2})_{\mathbb{Q}}-6h_{1}h_{2},$$and by point 2 of Theorem \ref{thm:m10} we have
$$B_{10}(f_{v}(\alpha_{1},h_{1}\sigma),f_{v}(\alpha_{2},h_{2}\sigma))=B_{10}(\widetilde{\lambda}(\alpha_{1})+h_{1}c_{1}(\widetilde{\Sigma}),\widetilde{\lambda}(\alpha_{2})+h_{2}c_{1}(\widetilde{\Sigma}))=$$ $$=(\alpha_{1},\alpha_{2})_{\mathbb{Q}}-6h_{1}h_{2},$$where by abuse of notation we denote $b_{v}$ and $B_{10}$ the $\mathbb{Q}-$bilinear extensions of the symmetric bilinear form $b_{v}$ on $(v^{\perp})^{*}\oplus_{\perp}\mathbb{Z}\cdot(\sigma/2)$ and of the Beauville form $B_{10}$ on $H^{2}(\widetilde{M}_{10},\mathbb{Z})$. It follows that the morphism $f_{v}$ sends the $\mathbb{Q}-$bilinear extension of $b_{v}$ to the $\mathbb{Q}-$bilinear extension of $B_{10}$. 

The morphism $f_{v}$ is an isomorphism of $\mathbb{Q}-$vector spaces by Theorem \ref{thm:m10}: it follows that $H^{2}(\widetilde{M}_{10},\mathbb{Z})$ is isometric to $f_{v}^{-1}(H^{2}(\widetilde{M}_{10},\mathbb{Z}))$. Finally, we have $f_{v}^{-1}(H^{2}(\widetilde{M}_{10},\mathbb{Z}))=\Gamma_{v}$. Indeed, let $(\alpha,k\sigma/2)\in\Gamma_{v}$: as $v=(2,0,-2)$ and $\alpha\in (v^{\perp})^{*}$, we have $\alpha=(n/2,\alpha',n/2)$ for some $n\in\mathbb{Z}$ and $\alpha'\in H^{2}(S,\mathbb{Z})$, hence $n\in 2\mathbb{Z}$ if and only if $k\in 2\mathbb{Z}$, i. e. $n+k$ is always even. By point 3 of Theorem \ref{thm:m10} we have $$f_{v}(\alpha,k\sigma/2)=\widetilde{\mu}(\alpha')+nc_{1}(\widetilde{B})+\frac{n+k}{2}c_{1}(\widetilde{\Sigma})\in H^{2}(\widetilde{M}_{10},\mathbb{Z}),$$so that $\Gamma_{v}\subseteq f_{v}^{-1}(H^{2}(\widetilde{M}_{10},\mathbb{Z})$. For the opposite inclusion, by point 1 of Theorem \ref{thm:m10} an element $\gamma\in H^{2}(\widetilde{M}_{10},\mathbb{Z})$ is of the form $\widetilde{\mu}(\alpha')+nc_{1}(\widetilde{B})+mc_{1}(\widetilde{\Sigma})$ for some $\alpha'\in H^{2}(S,\mathbb{Z})$ and $n,m\in\mathbb{Z}$. Hence we have $\gamma=f_{v}((n/2,\alpha',n/2),(m-n/2)\sigma)$, and $((n/2,\alpha',n/2),(m-n/2)\sigma)\in\Gamma_{v}$.

In particular, the symmetric bilinear form on $\Gamma_{v}$ is non-degenerate and takes values in $\mathbb{Z}$.}
\end{oss}

\begin{oss}
\label{oss:abe}{\rm If $S$ is abelian and $v=(2,0,-2)$, then it is easy to describe the lattice $H^{2}(\widetilde{K}_{6},\mathbb{Z})$ in terms of $v$: we consider the $\mathbb{Z}-$module $v^{\perp}\oplus_{\perp}\mathbb{Z}\cdot\alpha$, with the bilinear symmetric form $b_{v}$ defined by letting $b_{v|v^{\perp}}$ to be the Mukai pairing on $v^{\perp}$ and $b_{v}(\alpha,\alpha)=-2$. Then we have a morphism of $\mathbb{Z}-$modules $$f_{v}:v^{\perp}\oplus_{\perp}\mathbb{Z}\cdot\alpha\longrightarrow H^{2}(\widetilde{K}_{6},\mathbb{Z}),\,\,\,\,\,\,\,\,\,\,\,f_{v}(\beta,k\alpha):=\widetilde{\nu}(\beta)+kc_{1}(A).$$Using Theorem \ref{thm:k6} it is easy to show that $f_{v}$ is an isometry.}
\end{oss}

We now state and prove the main result of this section, in which we describe the pure weight-two Hodge structure and the lattice structure of $H^{2}(\widetilde{M}_{v},\mathbb{Z})$ (if $S$ is K3) and $H^{2}(\widetilde{K}_{v},\mathbb{Z})$ for every OLS-triple $(S,v,H)$. Before doing this, we recall that if $(S,v,H)$ is an OLS-triple, then $v=2w$ for some primitive Mukai vector $w$ such that $w^{2}=2$. On $v^{\perp}$ we have a pure weight-two Hodge structure $$(v^{\perp})^{0,2}:=(v^{\perp}\otimes\mathbb{C})\cap\widetilde{H}^{0,2}(S),\,\,\,\,\,\,\,\,\,\,(v^{\perp})^{2,0}:=(v^{\perp}\otimes\mathbb{C})\cap\widetilde{H}^{2,0}(S),$$ $$(v^{\perp})^{1,1}:=(v^{\perp}\otimes\mathbb{C})\cap\widetilde{H}^{1,1}(S).$$Notice that even the dual lattice $(v^{\perp})^{*}$ has then a pure weight-two Hodge structure, as $(v^{\perp})^{*}$ is a $\mathbb{Z}-$submodule of maximal rank of $v^{\perp}\otimes\mathbb{Q}$. Finally, we have a pure weight-two Hodge structure on $\Gamma_{v}$ as follows:
$$\Gamma_{v}^{0,2}:=(v^{\perp})^{0,2},\,\,\,\,\,\,\,\,\,\,\,\Gamma_{v}^{2,0}:=(v^{\perp})^{2,0},$$ $$\Gamma_{v}^{1,1}:=(v^{\perp})^{1,1}\oplus\mathbb{C}\cdot\sigma.$$Moreover, consider the $\mathbb{Z}-$module $v^{\perp}\oplus_{\perp}\mathbb{Z}\cdot A$ with the symmetric bilinear form $b_{v}$ defined by letting $b_{v|v^{\perp}}$ to be the Mukai pairing on $v^{\perp}$ and $b_{v}(A,A):=-2$. This is a lattice which carries a natural pure weight-two Hodge structure as follows:
$$(v^{\perp}\oplus_{\perp}\mathbb{Z}\cdot A)^{0,2}:=(v^{\perp})^{0,2},\,\,\,\,\,\,\,\,\,\,\,\,(v^{\perp}\oplus_{\perp}\mathbb{Z}\cdot A)^{2,0}:=(v^{\perp})^{2,0},$$ $$(v^{\perp}\oplus_{\perp}\mathbb{Z}\cdot A)^{1,1}:=(v^{\perp})^{1,1}\oplus\mathbb{C}\cdot A.$$The main result of this section is the following:

\begin{thm}
\label{thm:w2hs}Let $(S,v,H)$ be an OLS-triple.
\begin{enumerate}
 \item If $S$ is K3, then there is a Hodge isometry $$f_{v}:\Gamma_{v}\longrightarrow H^{2}(\widetilde{M}_{v},\mathbb{Z}).$$
 \item If $S$ is abelian, then there is a Hodge isometry $$f_{v}:v^{\perp}\oplus_{\perp}\mathbb{Z}\cdot A\longrightarrow H^{2}(\widetilde{K}_{v},\mathbb{Z}).$$
\end{enumerate}
\end{thm}

\proof We use the following notations: first of all, we write $v=2(r,\xi,a)$. Moreover, if $f:Y\longrightarrow T$ is any morphism of schemes and $L\in Pic(Y)$, then we write $Y_{t}:=f^{-1}(t)$ and $L_{t}:=L_{|Y_{t}}\in Pic(Y_{t})$. We divide the proof in two cases:

\textit{Case 1:} $S$ \textit{is a projective K3 surface}. In this case we have a natural morphism $$f_{v}:\Gamma_{v}\longrightarrow H^{2}(\widetilde{M}_{v},\mathbb{Q}),\,\,\,\,\,\,\,\,\,f_{v}(\alpha,k\sigma/2):=\widetilde{\lambda}_{v}(\alpha)+\frac{k}{2}c_{1}(\widetilde{\Sigma}_{v}),$$where here, by abuse of notation, we still denote $\widetilde{\lambda}_{v}$ the $\mathbb{Q}-$linear extension of $\widetilde{\lambda}_{v}$ (recall that $\alpha\in(v^{\perp})^{*}\subseteq v^{\perp}\otimes\mathbb{Q}$). Notice that $f_{v}$ is a Hodge morphism, as $\sigma$ and $c_{1}(\widetilde{\Sigma}_{v})$ are both $(1,1)-$classes and $\widetilde{\lambda}_{v}$ is a Hodge morphism. We are then left with verifying the following properties:
\begin{enumerate}
 \item $f_{v}:\Gamma_{v}\longrightarrow H^{2}(\widetilde{M}_{v},\mathbb{Z})$;
 \item $f_{v}$ is an isomorphism of $\mathbb{Z}-$modules;
 \item $f_{v}$ is an isometry.
\end{enumerate}
First of all, all these properties are verified in the case of $\widetilde{M}_{10}$, as shown in Remark \ref{oss:isom10}. Now, recall that $\widetilde{M}_{v}$ is deformation equivalent to $\widetilde{M}_{10}$, and this equivalence is obtained by using two kinds of transformations (see the proof of Theorem 1.6 of \cite{PR}): one is the deformation of $\widetilde{M}_{v}$ induced by a deformation of the OLS-triple $(S,v,H)$ along a smooth, connected curve $T$; the other is the isomorphism induced by some Fourier-Mukai transforms. It is then sufficient to show that the previous properties are stable under these transformations.

For the first kind of transformations, consider $T$ to be a smooth, connected curve, and $(\mathscr{X},\mathscr{L},\mathscr{H})$ a deformation of the OLS-triple $(S,v,H)$ along $T$. This is given by a smooth, projective family $\varphi:\mathscr{X}\longrightarrow T$ of K3 surfaces such that $\mathscr{X}_{0}\simeq S$ for some $0\in T$, and by two line bundles $\mathscr{L},\mathscr{H}\in Pic(\mathscr{X})$ such that $\mathscr{H}_{0}\simeq H$ and $\mathscr{L}_{0}\simeq L$, where $L\in Pic(S)$ is such that $c_{1}(L)=\xi$. For every $t\in T$ we write $\xi_{t}:=c_{1}(\mathscr{L}_{t})$, $v_{t}:=2(r,\xi_{t},a)$, and $T'$ for the subset of $T$ given by all the $t\in T$ such that $(\mathscr{X}_{t},v_{t},\mathscr{H}_{t})$ is an OLS-triple. 

We write $\phi:\mathscr{M}\longrightarrow T$ (resp. $\phi^{s}:\mathscr{M}^{s}\longrightarrow T$) for the relative moduli space of semistable (resp. stable) sheaves, $\Sigma:=\mathscr{M}\setminus\mathscr{M}^{s}$ and $\psi:\widetilde{\mathscr{M}}\longrightarrow T$ for the blow-up of $\mathscr{M}$ along $\Sigma$ with reduced structure. By \cite{PR} we know that $\widetilde{\mathscr{M}}_{t}=\widetilde{M}_{v_{t}}(\mathscr{X}_{t},\mathscr{H}_{t})$ for every $t\in T'$, and $R^{2}\psi_{*}\mathbb{Q}$ is a local system on $T$, such that $(R^{2}\psi_{*}\mathbb{Q})_{t}=H^{2}(\widetilde{M}_{v_{t}},\mathbb{Q})$ for every $t\in T'$.

Moreover, we have a local system $V$ on $T$ such that $V_{t}\simeq v_{t}^{\perp}$ for every $t\in T$. This allows us to define a local system $\Gamma\subseteq (V\otimes\mathbb{Q})\oplus_{\perp}\mathbb{Q}\cdot(\sigma/2)$ on $T$ such that $\Gamma_{t}=\Gamma_{v_{t}}$ for every $t\in T$. By the construction of the morphism $\widetilde{\lambda}_{v}$ given in \cite{PR}, we have a morphism $$\widetilde{\lambda}:V\longrightarrow R^{2}\psi_{*}\mathbb{Q},$$such that $\widetilde{\lambda}_{t}=\widetilde{\lambda}_{v_{t}}$ for every $t\in T'$. This allows us to define a morphism $$f:\Gamma\longrightarrow R^{2}\psi_{*}\mathbb{Q}$$such that $f_{t}=f_{v_{t}}$ for every $t\in T'$. Hence $f_{v}$ takes values in $H^{2}(\widetilde{M}_{v},\mathbb{Z})$ (resp. it is an isomorphism of $\mathbb{Z}-$modules, it is an isometry) if and only if there is $t\in T'$ such that $f_{v_{t}}$ does.

For the second kind of transformation, let $(S,v,H)$ and $(S,v',H')$ be two OLS-triples, and suppose there is a Fourier-Mukai transform $FM:D^{b}(S)\longrightarrow D^{b}(S)$ verifying the two following properties:
\begin{enumerate}
 \item for every $H-$stable sheaf $\mathscr{E}$ of Mukai vector $v$ (resp. $w:=v/2$) we have that $FM(\mathscr{E})$ is $H'-$stable and has Mukai vector $v'$ (resp. $w':=v'/2$);
 \item for every $H'-$stable sheaf $\mathscr{E}$ of Mukai vector $v'$ (resp. $w'$) we have that $FM^{-1}(\mathscr{E})$ is $H-$stable and has Mukai vector $v$ (resp. $w$).
\end{enumerate}
Then $FM$ induces an isomorphism $g:M_{v}\longrightarrow M_{v'}$ which gives an isomorphism between $\Sigma_{v}$ and $\Sigma_{v'}$. Hence, it induces an isomorphism $\widetilde{g}:\widetilde{M}_{v}\longrightarrow\widetilde{M}_{v'}$, and $$\widetilde{g}^{*}:H^{2}(\widetilde{M}_{v'},\mathbb{Z})\longrightarrow H^{2}(\widetilde{M}_{v},\mathbb{Z})$$is a Hodge isometry such that $\widetilde{g}^{*}(c_{1}(\widetilde{\Sigma}_{v'}))=c_{1}(\widetilde{\Sigma}_{v})$. Moreover, $FM$ induces a Hodge isometry $h:(v')^{\perp}\longrightarrow v^{\perp}$, and the following diagram
\begin{equation}
\label{eq:comm}
\begin{array}{ccc}
(v')^{\perp} & \stackrel{h}\longrightarrow & v^{\perp}\\
\scriptstyle{\lambda_{v'}}\downarrow & & \downarrow\scriptstyle{\lambda_{v}}\\
H^{2}(M_{v'},\mathbb{Z}) & \stackrel{g^{*}}\longrightarrow & H^{2}(M_{v},\mathbb{Z})\\
\scriptstyle{\pi_{v'}^{*}}\downarrow & & \downarrow\scriptstyle{\pi_{v}^{*}}\\
H^{2}(\widetilde{M}_{v'},\mathbb{Z}) & \stackrel{\widetilde{g}^{*}}\longrightarrow & H^{2}(\widetilde{M}_{v},\mathbb{Z})
\end{array}
\end{equation}
is commutative: the commutativity of the top diagram is shown in Lemma 3.11 of \cite{PR}, and the commutativity of the down diagram comes from the definition of $\widetilde{g}$. Notice that the Hodge isometry $h$ extends to a Hodge isometry $$\widetilde{h}:\Gamma_{v'}\longrightarrow\Gamma_{v},\,\,\,\,\,\,\,\,\widetilde{h}(\alpha,k\sigma/2):=(h(\alpha),k\sigma/2),$$where, by abuse of notation, we still write $h$ for the $\mathbb{Q}-$linear extension of $h$. The commutativity of the diagram (\ref{eq:comm}) then implies that the following diagram
$$\begin{array}{ccc}
\Gamma_{v'} & \stackrel{\widetilde{h}}\longrightarrow & \Gamma_{v}\\
\scriptstyle{f_{v'}}\downarrow & & \downarrow\scriptstyle{f_{v}}\\
H^{2}(\widetilde{M}_{v'},\mathbb{Q}) & \stackrel{\widetilde{g}^{*}}\longrightarrow & H^{2}(\widetilde{M}_{v},\mathbb{Q})
\end{array}$$is commutative. As $\widetilde{h}$ and $\widetilde{g}^{*}$ are two isometries and $\widetilde{g}$ is an isomorphism, then $f_{v}$ takes values in $H^{2}(\widetilde{M}_{v},\mathbb{Z})$ (resp. it is an isomorphism of $\mathbb{Z}-$modules, it is an isometry) if and only if $f_{v'}$ does.

\textit{Case 2:} $S$ \textit{is an abelian surface}. We have a natural morphism $$f_{v}:v^{\perp}\oplus_{\perp}\mathbb{Z}\cdot A\longrightarrow H^{2}(\widetilde{K}_{v},\mathbb{Q}),\,\,\,\,\,\,\,\,\,\,f_{v}(\alpha,kA):=\widetilde{\nu}_{v}(\alpha)+\frac{k}{2}c_{1}(\widetilde{\Sigma}_{v}),$$which is easily seen to be a Hodge morphism. We then just need to show that $f_{v}$ takes values in $H^{2}(\widetilde{K}_{v},\mathbb{Z})$, that is an isomorphism of $\mathbb{Z}-$modules, and that it is an isometry. By Remark \ref{oss:abe} these properties are verified in the case of $\widetilde{K}_{6}$. To show them for every OLS-triple, the proof is formally the same as in the previous case.\endproof

\begin{oss}
\label{oss:divkv}{\rm A difference between $\widetilde{M}_{v}$ for a K3 surface and $\widetilde{K}_{v}$ for an abelian surface is that in the first case the exceptional divisor is primitive, while in the second case it is divisible by 2. This follows from Theorem \ref{thm:w2hs}. 

To see it more directly, recall that $\widetilde{M}_{v}$ (resp. $\widetilde{K}_{v}$) is deformation equivalent to $\widetilde{M}_{10}$ (resp. to $\widetilde{K}_{6}$), and the deformation equivalence is realized as recalled in the proof of Theorem \ref{thm:w2hs}. Following this deformation, we notice that the exceptional divisor $\widetilde{\Sigma}_{v}$ of $\widetilde{M}_{v}$ (resp. of $\widetilde{K}_{v}$) deforms to the exceptional divisor $\widetilde{\Sigma}$ of $\widetilde{M}_{10}$ (resp. $\widetilde{K}_{6}$). The divisibility of $\widetilde{\Sigma}_{v}$ is then equivalent to the divisibility of $\widetilde{\Sigma}$: on $\widetilde{M}_{10}$, we have that $\widetilde{\Sigma}$ is primitive by Theorem \ref{thm:m10}, while on $\widetilde{K}_{6}$ the divisor $\widetilde{\Sigma}$ is divisible by 2 by Theorem \ref{thm:k6}. The square root of $\mathscr{O}(\widetilde{\Sigma}_{v})$ plays in important role for the $2-$factoriality of $K_{v}$, and we denote it $A_{v}$.}
\end{oss}

\section{Factoriality of $M_{v}(S,H)$ and $K_{v}(S,H)$}

The aim of this section is to prove Theorems \ref{thm:maink3} and \ref{thm:mainab}. We use the following notation: if $L$ is a $\mathbb{Z}-$module with a pure Hodge structure of weight 2, we write $L^{1,1}_{\mathbb{Z}}:=L\cap L^{1,1}$.

\subsection{The proof of Theorem \ref{thm:maink3}}

We begin with the factoriality of the moduli space $M_{v}(S,H)$ associated to an OLS-triple $(S,v,H)$ where $S$ is a K3 surface. First of all, we prove the following result in which we present a criterion for the local factoriality and the $2-$factoriality of $M_{v}(S,H)$. 

\begin{prop}
\label{prop:crit1} 
Let $(S,v,H)$ be an OLS-triple such that $S$ is a projective K3 surface. Then $M_{v}(S,H)$ is either $2-$factorial or locally factorial. Moreover
 \begin{enumerate}
  \item $M_{v}(S,H)$ is $2-$factorial if and only if there is $\beta\in(v^{\perp})^{1,1}_{\mathbb{Z}}$ primitive such that $(\beta,v^{\perp})\subseteq 2\mathbb{Z}$.
  \item $M_{v}(S,H)$ is locally factorial if and only if for every primitive $\beta\in(v^{\perp})^{1,1}_{\mathbb{Z}}$ we have $(\beta,v^{\perp})=\mathbb{Z}$.
 \end{enumerate}
\end{prop}

\proof We know that $\widetilde{M}_{v}$ is simply connected, so that $Pic(\widetilde{M}_{v})=H^{2}(\widetilde{M}_{v},\mathbb{Z})\cap H^{1,1}(\widetilde{M}_{v})$. By Remark \ref{oss:explicitgamma} we have that $(\Gamma_{v})^{1,1}_{\mathbb{Z}}$ is equal to the set $$\{(\beta/2,k\sigma/2)\,|\,\beta\in (v^{\perp})^{1,1}_{\mathbb{Z}},\,\,(\beta,v^{\perp})\subseteq 2\mathbb{Z},\,\,k\in\mathbb{Z},\,\,\mathrm{and}\,\,k\in 2\mathbb{Z}\,\,\mathrm{iff}\,\,\beta/2\in v^{\perp}\}.$$By point 1 of Theorem \ref{thm:w2hs} we have $Pic(\widetilde{M}_{v})=f_{v}((\Gamma_{v})^{1,1}_{\mathbb{Z}})$.

We have an exact sequence $$0\longrightarrow\mathbb{Z}\stackrel{\iota}\longrightarrow Pic(\widetilde{M}_{v})\stackrel{r}\longrightarrow Pic(M_{v}^{s})\longrightarrow 0,$$where $\iota(1):=\mathscr{O}(\widetilde{\Sigma}_{v})$ and $r$ is the restriction morphism. From this exact sequence it follows that $$Pic(M_{v}^{s})=Pic(\widetilde{M}_{v})/\mathbb{Z}\cdot\mathscr{O}(\widetilde{\Sigma}_{v})=f_{v}((\Gamma_{v})^{1,1}_{\mathbb{Z}})/\mathbb{Z}\cdot f_{v}(\sigma).$$Since the complement of $M_{v}^{s}$ in $M_{v}$ has codimension 2, we have $A^{1}(M_{v})=Pic(M_{v}^{s})$. Moreover, by Corollary \ref{cor:picmv} we have $Pic(M_{v})=L_{v}((v^{\perp})^{1,1}_{\mathbb{Z}})$, hence $\pi_{v}^{*}(Pic(M_{v}))=f_{v}((v^{\perp})^{1,1}_{\mathbb{Z}})$. Therefore $$A^{1}(M_{v})/Pic(M_{v})= Pic(\widetilde{M}_{v})/(f_{v}((v^{\perp})^{1,1}_{\mathbb{Z}})\oplus_{\perp}\mathbb{Z}\cdot f_{v}(\sigma))\simeq$$ $$\simeq (\Gamma_{v})^{1,1}_{\mathbb{Z}}/((v^{\perp})^{1,1}_{\mathbb{Z}}\oplus_{\perp}\mathbb{Z}\cdot\sigma).$$This quotient is either $\mathbb{Z}/2\mathbb{Z}$ or trivial. It is $\mathbb{Z}/2\mathbb{Z}$, i. e. $M_{v}$ is $2-$factorial, if and only if there is $\beta'\in(v^{\perp})^{1,1}_{\mathbb{Z}}$ which is not divisible by 2 in $(v^{\perp})^{1,1}_{\mathbb{Z}}$ such that $(\beta',v^{\perp})\subseteq 2\mathbb{Z}$. If $N:=max\{n\in\mathbb{N}\,|\,\beta'/n\in(v^{\perp})^{1,1}_{\mathbb{Z}}\}$, then $N$ is odd and $\beta:=\beta'/N$ is a primitive element of $(v^{\perp})^{1,1}_{\mathbb{Z}}$ such that $(\beta,v^{\perp})\subseteq 2\mathbb{Z}$.

The quotient $(\Gamma_{v})^{1,1}_{\mathbb{Z}}/((v^{\perp})^{1,1}_{\mathbb{Z}}\oplus_{\perp}\mathbb{Z}\cdot\sigma)$ is trivial, i. e. $M_{v}$ is locally factorial, if and only if for every $\beta\in(v^{\perp})^{1,1}_{\mathbb{Z}}$ such that $(\beta,v^{\perp})\subseteq 2\mathbb{Z}$ there is $\alpha\in(v^{\perp})^{1,1}_{\mathbb{Z}}$ such that $\beta=2\alpha$. This means that for every primitive $\beta\in(v^{\perp})^{1,1}_{\mathbb{Z}}$ we cannot have $(\beta,v^{\perp})\subseteq 2\mathbb{Z}$: as $|\Delta(v^{\perp})|=2$, this implies $(\beta,v^{\perp})=\mathbb{Z}$.\endproof

\begin{oss}
\label{oss:ugualos}
{\rm Notice that as $|\Delta(v^{\perp})|=2$, it follows that if $\beta\in (v^{\perp})^{1,1}_{\mathbb{Z}}$ is such that $(\beta,v^{\perp})\subseteq 2\mathbb{Z}$, then we even have $(\beta,v^{\perp})=2\mathbb{Z}$.}
\end{oss}

The statement of Proposition \ref{prop:crit1} has an equivalent reformulation which gives a criterion for the local factoriality and for the 2-factoriality which is easier to verify in concrete situations.
\par\bigskip
\noindent\textbf{Theorem \ref{thm:maink3}}. \textit{Let} $(S,v,H)$ \textit{be an OLS-triple such that} $S$ \textit{is a projective K3 surface, and write} $v=2w$\textit{. Then} $M_{v}(S,H)$ \textit{is either} $2-$\textit{factorial or locally factorial. Moreover}
 \begin{enumerate}
  \item $M_{v}(S,H)$ \textit{is} $2-$\textit{factorial if and only if there is} $\gamma\in \widetilde{H}(S,\mathbb{Z})\cap\widetilde{H}^{1,1}(S)$ \textit{such that} $(\gamma,w)=1$\textit{.}
  \item $M_{v}(S,H)$ \textit{is locally factorial if and only if for every} $\gamma\in\widetilde{H}(S,\mathbb{Z})\cap\widetilde{H}^{1,1}(S)$ \textit{we have} $(\gamma,w)\in 2\mathbb{Z}$\textit{.}
 \end{enumerate}

\proof The fact that $M_{v}$ is either locally factorial of $2-$factorial was already shown in Proposition \ref{prop:crit1}. We are then left to show the two criteria for the $2-$factoriality and for the local factoriality. As $w^{2}=2$, then either there is $\gamma\in\widetilde{H}^{1,1}(S)\cap H^{2}(S,\mathbb{Z})$ such that $(\gamma,w)=1$, or for every $\gamma\in\widetilde{H}(S)\cap H^{2}(S,\mathbb{Z})$ we have $(\gamma,w)\in 2\mathbb{Z}$. Hence, we just need to show that $M_{v}$ is $2-$factorial if and only if there is $\gamma\in\widetilde{H}(S,\mathbb{Z})\cap\widetilde{H}^{1,1}(S)$ such that $(\gamma,w)=1$.

As $w^{2}=2$, we have that $|\Delta(v^{\perp})|=|\Delta(\mathbb{Z}\cdot w)|=2$. It follows that $v^{\perp}\oplus_{\perp}\mathbb{Z}\cdot w$ is a sublattice of $\widetilde{H}(S,\mathbb{Z})$ of index 2. Moreover, the natural inclusion of $\widetilde{H}(S,\mathbb{Z})$ in $(v^{\perp})^{*}\oplus_{\perp}(\mathbb{Z}\cdot w)^{*}$ presents $\widetilde{H}(S,\mathbb{Z})$ as a sublattice of index 2 of $(v^{\perp})^{*}\oplus_{\perp}(\mathbb{Z}\cdot w)^{*}$. More precisely, if $$H_{v}:=\{(\delta_{1},\delta_{2})\in (v^{\perp})^{*}\oplus_{\perp}(\mathbb{Z}\cdot w)^{*}\,|\,\delta_{1}\in v^{\perp}\,\,\mathrm{if}\,\,\mathrm{and}\,\,\mathrm{only}\,\,\mathrm{if}\,\,\delta_{2}\in\mathbb{Z}\cdot w\},$$then $\widetilde{H}(S,\mathbb{Z})=H_{v}$ (indeed, we have $((v^{\perp})^{*}\oplus_{\perp}(\mathbb{Z}\cdot w)^{*})/(v^{\perp}\oplus_{\perp}\mathbb{Z}\cdot w)\simeq(\mathbb{Z}/2\mathbb{Z})^{2}$, hence there are only three sublattices of $(v^{\perp})^{*}\oplus_{\perp}(\mathbb{Z}\cdot w)^{*}$ of index 2, which are $(v^{\perp})^{*}\oplus_{\perp}\mathbb{Z}\cdot w$, $v^{\perp}\oplus_{\perp}(\mathbb{Z}\cdot w)^{*}$ and $H_{v}$; but notice that $((v^{\perp})^{*}\oplus_{\perp}\mathbb{Z}\cdot w)\cap\widetilde{H}(S,\mathbb{Z})$ and $(v^{\perp}\oplus_{\perp}(\mathbb{Z}\cdot w)^{*})\cap\widetilde{H}(S,\mathbb{Z})$ are both equal to $v^{\perp}\oplus_{\perp}\mathbb{Z}\cdot w$, hence we have $\widetilde{H}(S,\mathbb{Z})=H_{v}$).

Now, let us suppose that $M_{v}$ is $2-$factorial. By Proposition \ref{prop:crit1} it follows that there is a primitive $\beta\in(v^{\perp})^{1,1}_{\mathbb{Z}}$ such that $(\beta,v^{\perp})\subseteq 2\mathbb{Z}$. This implies that $\beta/2\in (v^{\perp})^{*}$: then $\gamma:=(\beta/2,w/2)\in H_{v}$, i. e. $\gamma\in\widetilde{H}(S,\mathbb{Z})$. As both $\beta$ and $w$ are $(1,1)-$classes, it follows that $\gamma\in\widetilde{H}(S,\mathbb{Z})\cap\widetilde{H}^{1,1}(S)$. Moreover, as $\beta\in v^{\perp}$ and $w^{2}=2$ we have $(\gamma,w)=w^{2}/2=1.$

Conversely, suppose that there is $\gamma\in\widetilde{H}(S,\mathbb{Z})\cap\widetilde{H}^{1,1}(S)$ such that $(\gamma,w)=1$. By the equality $\widetilde{H}(S,\mathbb{Z})=H_{v}$ there are $\delta_{1}\in(v^{\perp})^{*}$ and $\delta_{2}\in(\mathbb{Z}\cdot w)^{*}$ such that $\gamma=(\delta_{1},\delta_{2})$. As $\gamma$ is a $(1,1)-$class, we even have that $\delta_{1}$ has to be a $(1,1)-$class in $(v^{\perp})^{*}$. As $(\delta_{1},w)=0$, it follows that $\delta_{2}=w/2$, so that $\delta_{1}\notin v^{\perp}$. Now, recall that $v^{\perp}$ has index 2 in $(v^{\perp})^{*}$: it follows that $\beta':=2\delta_{1}\in (v^{\perp})^{1,1}_{\mathbb{Z}}$ is not divisible by 2 in $v^{\perp}$, and $(\beta',v^{\perp})\subseteq 2\mathbb{Z}$. If $N:=max\{n\in\mathbb{N}\,|\,\beta'/N\in v^{\perp}\}$, then $N$ is odd and $\beta:=\beta'/N\in(v^{\perp})^{1,1}_{\mathbb{Z}}$ is primitive, and $(\beta,v^{\perp})\subseteq 2\mathbb{Z}$. By Proposition \ref{prop:crit1} it follows that $M_{v}$ is $2-$factorial.
\endproof

\begin{oss}
\label{oss:muk}
{\rm We observe that the condition which guarantees that $M_{v}(S,H)$ is $2-$factorial in Theorem \ref{thm:maink3} is the same one that Mukai shows to imply the existence of a universal family on $M_{w}(S,H)$ (see \cite{M1}).}
\end{oss}

From Theorem \ref{thm:maink3} it is easy to produce examples of OLS-triples $(S,v,H)$ such that $M_{v}(S,H)$ is locally factorial:

\begin{esem}
\label{esem:lfmv}
{\rm Let $n,m\in\mathbb{Z}$, and $S$ be a K3 surface with $Pic(S)=\mathbb{Z}$, such that there is $\xi\in Pic(S)$ with $\xi^{2}=2+2nm$. Then the Mukai vector $w:=(n,\xi,m)$ is such that $w^{2}=2$, hence it is primitive. Moreover, as $Pic(S)=\mathbb{Z}$, then for every $\zeta\in Pic(S)$ we have $\xi\cdot\zeta\in 2\mathbb{Z}$: it follows that $M_{2w}$ is locally factorial if $n,m\in 2\mathbb{Z}$. 

If one between $m$ and $n$ is odd, then $M_{2w}$ is $2-$factorial. Indeed, suppose that $m$ is odd: as $w^{2}=2$ then $m$ has to be prime with $\xi^{2}$. It follows that there are $h,k\in\mathbb{Z}$ such that $((0,h\xi,k),w)=1$, and Theorem \ref{thm:maink3} implies that $M_{2w}$ is $2-$factorial.}
\end{esem}

\subsection{The proof of Theorem \ref{thm:mainab}}

This section is devoted to the study of the factoriality properties of $M_{v}(S,H)$ and $K_{v}(S,H)$ for every OLS-triple $(S,v,H)$ where $S$ is abelian. More precisely, we prove
\par\bigskip
\noindent\textbf{Theorem 1.2.} \textit{Let} $(S,v,H)$ \textit{be an OLS-triple such that} $S$ \textit{is an abelian surface. Then} $M_{v}(S,H)$ \textit{and} $K_{v}(S,H)$ \textit{are} $2-$\textit{factorial.}
\par\bigskip
Before starting the proof, we introduce here some notations which will be useful in the following, namely the commutative diagrams:
\begin{equation}
\label{eq:dued}
\begin{array}{ccc}
K_{v}^{s} & \stackrel{i_{v}^{0}}\longrightarrow & K_{v}\\
\scriptstyle{j_{v}^{s}}\downarrow & & \downarrow\scriptstyle{j_{v}}\\
M_{v}^{s} & \stackrel{i_{v}}\longrightarrow & M_{v}
\end{array}\,\,\,\,\,\,\,\,\,\,\mathrm{and}\,\,\,\,\,\,\,\,\,
\begin{array}{ccc}
K_{v}^{s} & \stackrel{\widetilde{i}_{v}^{0}}\longrightarrow & \widetilde{K}_{v}\\
\scriptstyle{j_{v}^{s}}\downarrow & & \downarrow\scriptstyle{\widetilde{j}_{v}}\\
M_{v}^{s} & \stackrel{\widetilde{i}_{v}}\longrightarrow & \widetilde{M}_{v}
\end{array}
\end{equation}
where $i_{v},i_{v}^{0},j_{v},j_{v}^{s},\widetilde{j}_{v}$ are the inclusions, and where $\widetilde{i}_{v}$ (resp. $\widetilde{i}_{v}^{0}$) is the composition of $(\pi_{v}^{-1})_{|M_{v}^{s}}$ (resp. $(\pi_{v}^{-1})_{|K_{v}^{s}}$) with the inclusion of $\pi_{v}^{-1}(M_{v}^{s})$ (resp. of $\pi_{v}^{-1}(K_{v}^{s})$) in $\widetilde{M}_{v}$ (resp. $\widetilde{K}_{v}$). Moreover, we let $\widetilde{E}$ be the exceptional divisor of $\pi_{v}:\widetilde{M}_{v}\longrightarrow M_{v}$. Notice that $\widetilde{j}_{v}^{*}(\widetilde{E})=\widetilde{\Sigma}_{v}$.

\subsubsection{The $2-$factoriality of $K_{v}(S,H)$}

We begin by proving that $K_{v}(S,H)$ is $2-$factorial, showing then the first half of Theorem \ref{thm:mainab}. The proof of the $2-$factoriality of $K_{v}(S,H)$ is similar to the one of Proposition \ref{prop:crit1}. Moreover, here we want to stress the role of the line bundle $A_{v}\in Pic(\widetilde{K}_{v})$ such that $A_{v}^{\otimes 2}=\mathscr{O}(\widetilde{\Sigma}_{v})$ (see Remark \ref{oss:divkv}). 

\begin{prop}
\label{prop:2f}Let $(S,v,H)$ be an OLS-triple such that $S$ is an abelian surface. Then $K_{v}(S,H)$ is $2-$factorial, and $Pic(K^{s}_{v})/Pic(K_{v})$ is generated by the class of $(\widetilde{i}^{0}_{v})^{*}(A_{v})$. 
\end{prop}

\proof As $\widetilde{K}_{v}$ is simply connected, we have $Pic(\widetilde{K}_{v})=H^{2}(\widetilde{K}_{v},\mathbb{Z})\cap H^{1,1}(\widetilde{K}_{v})$. By point 2 of Theorem \ref{thm:w2hs} we have then $$Pic(\widetilde{K}_{v})=f_{v}((v^{\perp})^{1,1}_{\mathbb{Z}}\oplus_{\perp}\mathbb{Z}\cdot A).$$Now, we have an exact sequence $$0\longrightarrow\mathbb{Z}\stackrel{\iota}\longrightarrow Pic(\widetilde{K}_{v})\stackrel{r}\longrightarrow Pic(K_{v}^{s})\longrightarrow 0,$$where $\iota(1):=\mathscr{O}(\widetilde{\Sigma}_{v})$ and $r$ is the restriction morphism. From this exact sequence and the fact that $c_{1}(\mathscr{O}(\widetilde{\Sigma}_{v}))=f_{v}(2A)$, it follows that $$Pic(K_{v}^{s})=Pic(\widetilde{K}_{v})/\mathbb{Z}\cdot\mathscr{O}(\widetilde{\Sigma}_{v})=f_{v}((v^{\perp})^{1,1}_{\mathbb{Z}}\oplus_{\perp}\mathbb{Z}\cdot A)/\mathbb{Z}\cdot f_{v}(2A).$$Since the complement of $K_{v}^{s}$ in $K_{v}$ has codimension 2, we have $A^{1}(K_{v})=Pic(K_{v}^{s})$. Moreover, by point 2 of Corollary \ref{cor:picmv} we have $Pic(K_{v})=N_{v}((v^{\perp})^{1,1}_{\mathbb{Z}})$, hence $\pi_{v}^{*}(Pic(K_{v}))=f_{v}((v^{\perp})^{1,1}_{\mathbb{Z}})$. Therefore $$A^{1}(K_{v})/Pic(K_{v})= Pic(\widetilde{K}_{v})/(f_{v}((v^{\perp})^{1,1}_{\mathbb{Z}})\oplus_{\perp}\mathbb{Z}\cdot f_{v}(2A))\simeq$$ $$\simeq((v^{\perp})^{1,1}_{\mathbb{Z}}\oplus_{\perp}\mathbb{Z}\cdot A)/((v^{\perp})^{1,1}_{\mathbb{Z}}\oplus_{\perp}\mathbb{Z}\cdot 2A)\simeq\mathbb{Z}/2\mathbb{Z},$$and the generator is the class of $(\widetilde{i}_{v}^{0})^{*}A_{v}$.\endproof

\subsubsection{The $2-$factoriality of $M_{v}(S,H)$}

By Proposition \ref{prop:2f}, to complete the proof of Theorem \ref{thm:mainab} it only remains to show that $M_{v}(S,H)$ is $2-$factorial for every OLS-triple $(S,v,H)$ where $S$ is abelian. First, by using Proposition \ref{prop:2f} we show that the moduli space $M_{v}(S,H)$ is either locally factorial or $2-$factorial.

\begin{prop}
\label{prop:mvlf2f}Let $(S,v,H)$ be an OLS-triple such that $S$ is an abelian surface. Then $M_{v}(S,H)$ is either locally factorial or $2-$factorial.
\end{prop}

\proof We have a finite étale cover $$K_{v}\times(S\times\widehat{S})\longrightarrow M_{v},\,\,\,\,\,\,\,\,\,(\mathscr{F},p,L)\mapsto t_{p}^{*}\mathscr{F}\otimes L,$$where $t_{p}:S\longrightarrow S$ is the translation by $p$. We use the following general fact in commutative algebra: if $A$ and $B$ are two Noetherian local rings and $f:A\longrightarrow B$ is a morphism of local rings giving to $B$ the structure of a flat $A-$module, then the natural morphism $Cl(A)\longrightarrow Cl(B)$ is injective (see Proposition 21.13.12 of EGA IV). Applying this to our situation, we see that if $K_{v}\times(S\times\widehat{S})$ is $2-$factorial, then $M_{v}$ is either locally factorial or $2-$factorial.

We are then left with showing that $K_{v}\times(S\times\widehat{S})$ is $2-$factorial. By Remark 5.3 of \cite{BGS} we have that $$A^{1}(K_{v}\times(S\times\widehat{S}))=A^{1}(K_{v})\times A^{1}(S\times\widehat{S})\times Hom(Alb(K_{v}),P(S\times\widehat{S})),$$where $Alb(K_{v})$ is the Albanese variety of the singular variety $K_{v}$, and $P(S\times\widehat{S})$ is the Picard variety of $S\times\widehat{S}$. Notice that as $\widetilde{K}_{v}$ is irreducible symplectic by point 2 of Theorem \ref{thm:pr}, we have that $Alb(\widetilde{K}_{v})$ is trivial, so that $Alb(K_{v})$ is trivial, and $$A^{1}(K_{v}\times(S\times\widehat{S}))=A^{1}(K_{v})\times A^{1}(S\times\widehat{S}).$$We have a similar description of the Picard group of $K_{v}\times(S\times\widehat{S})$ $$Pic(K_{v}\times(S\times\widehat{S}))=Pic(K_{v})\times Pic(S\times\widehat{S})\times Hom(Alb(K_{v}),P(S\times\widehat{S}))=$$ $$=Pic(K_{v})\times Pic(S\times\widehat{S}).$$As $S\times\widehat{S}$ is smooth, we have that $A^{1}(S\times\widehat{S})=Pic(S\times\widehat{S})$. It follows that $$A^{1}(K_{v}\times(S\times\widehat{S}))/Pic(K_{v}\times(S\times\widehat{S}))=A^{1}(K_{v})/Pic(K_{v})\simeq\mathbb{Z}/2\mathbb{Z}$$by Proposition \ref{prop:2f}. Hence $K_{v}\times(S\times\widehat{S})$ is $2-$factorial, and we are done.\endproof

To conclude the proof of Theorem \ref{thm:mainab}, by Proposition \ref{prop:mvlf2f} we have then just to show that $M_{v}$ is not locally factorial, i. e. that $Pic(M_{v})$ is strictly contained in $Pic(M_{v}^{s})$. This is true if there is a line bundle $L\in Pic(\widetilde{M}_{v})$ such that $\widetilde{j}_{v}^{*}(L)=A_{v}$. Indeed, consider $\widetilde{i}_{v}^{*}(L)\in Pic(M_{v}^{s})$: if there was $L'\in Pic(M_{v})$ such that $\widetilde{i}_{v}^{*}(L)=i_{v}^{*}(L')$, then by the commutativity of the diagrams (\ref{eq:dued}) we would have $$(\widetilde{i}_{v}^{0})^{*}A_{v}=(\widetilde{i}_{v}^{0})^{*}\widetilde{j}_{v}^{*}(L)=(i_{v}^{0})^{*}j_{v}^{*}(L').$$As $j_{v}^{*}(L')\in Pic(K_{v})$, it would follow that $(\widetilde{i}_{v}^{0})^{*}A_{v}$ is trivial in $Pic(K_{v}^{s})/Pic(K_{v})$, which is not possible from Proposition \ref{prop:2f}.

In conclusion, in order to complete the proof of Theorem \ref{thm:mainab} it is sufficient to find $L\in Pic(\widetilde{M}_{v})$ which restricts to $A_{v}$ on $\widetilde{K}_{v}$ or, equivalently, it is sufficient to find a line bundle on $Pic(\widetilde{M}_{v})$ whose restriction to $\widetilde{K}_{v}$ is $\mathscr{O}(\widetilde{\Sigma}_{v})$, and whose first Chern class is divisible by 2 in $H^{2}(\widetilde{M}_{v},\mathbb{Z})$.

The strategy to prove the existence of such a line bundle will be the following: first, we prove a criterion to determine if for a class $\alpha\in H^{2}(\widetilde{M}_{v},\mathbb{Z})$ there is $\beta\in H^{2}(\widetilde{M}_{v},\mathbb{Z})$ such that $\alpha=2\beta$; then we provide a line bundle on $\widetilde{M}_{v}$ whose first Chern class verifies such a criterion.

\paragraph{Criterion for the divisibility in $H^{2}(\widetilde{M}_{v},\mathbb{Z})$.}

In order to determine if a class $\alpha\in H^{2}(\widetilde{M}_{v},\mathbb{Z})$ is divisible by some integer $n$ in $H^{2}(\widetilde{M}_{v},\mathbb{Z})$, we will see that it is sufficient to show that $\alpha$ is divisible by $n$ when restricted to $\widetilde{K}_{v}$ and to a certain $\mathbb{P}^{1}-$bundle $\mathbb{P}\subseteq\widetilde{M}_{v}$. The first part of this section is devoted to the construction of such a $\mathbb{P}^{1}-$bundle.

Let $F\in M_{w}$, so that $F$ is an $H-$stable sheaf of $S$ with Mukai vector $w$. Moreover, we consider the following morphisms: we write $\pi_{1}$ and $\pi_{2}$ for the two natural projections of $S\times S$, and $p_{ij}$ for the projection of $S\times S\times\widehat{S}$ to the product of the $i-$th and of the $j-$th factors. Finally, consider the morphism $$f:S\times S\longrightarrow S\times S,\,\,\,\,\,\,\,\,\,\,\,f(p,q):=(p+q,q).$$Define $$\mathscr{F}:=p_{12}^{*}(f^{*}(\pi_{1}^{*}F))\otimes p_{13}^{*}\mathscr{P},$$where $\mathscr{P}$ is the Poincaré bundle on $S\times\widehat{S}$. It is easy to see that for every $p\in S$ and $L\in\widehat{S}$ we have that $$\mathscr{F}_{|p_{23}^{-1}(p,L)}\simeq t_{p}^{*}F\otimes L,$$where $$t_{p}:S\longrightarrow S,\,\,\,\,\,\,\,\,\,\,\,t_{p}(q):=q+p$$is the translation with $p$.

Now, recall that as $w$ is primitive and $w^{2}=2$, then by \cite{M2} and \cite{Y4} the morphism $a_{w}:M_{w}\longrightarrow S\times\widehat{S}$ is an isomorphism. Let $$\mathscr{U}:=(id_{S}\times a_{w})^{*}\mathscr{F},$$which is an $M_{w}-$flat family of coherent sheaves on $S\times M_{w}$. Moreover, it is a universal family: indeed, we have a morphism $$h:S\times\widehat{S}\longrightarrow M_{w},\,\,\,\,\,\,\,\,\,\,\,h(p,L):=t_{p}^{*}F\otimes L,$$as for every $p\in S$ and $L\in\widehat{S}$ the sheaf $t_{p}^{*}F\otimes L$ is $H-$stable of Mukai vector $w$. This is an isomorphism, as shown in section 1.2 of \cite{Y3}, and $\mathscr{U}$ is now easily seen to be a universal family on $S\times M_{w}$.

Now, consider $G\in M_{w}$ and let $\mathscr{U}':=\mathscr{U}_{|S\times(M_{w}\setminus\{G\})}$. Let $q_{1}$ and $q_{2}$ be the two projections of $S\times (M_{w}\setminus\{G\})$ respectively to $S$ and to $M_{w}\setminus\{G\}$, and define $$\mathscr{V}:=\mathscr{E}xt^{1}_{q_{2}}(q_{1}^{*}G,\mathscr{U}').$$This is a locally free sheaf of rank 2, as showed in the following:

\begin{lem}
\label{lem:vb2}The sheaf $\mathscr{V}$ is a vector bundle of rank 2 over $M_{w}\setminus\{G\}$, and for every $H\in M_{w}\setminus\{G\}$ we have $\mathscr{V}_{H}\simeq Ext^{1}(G,H)$.
\end{lem}

\proof As $G$ and $H$ are two non-isomorphic stable sheaves with Mukai vector $w$, we have that $\chi(F,G)=-2$ and $Hom(G,H)=Ext^{2}(G,H)=0$. This implies that $Ext^{1}(G,H)$ is a $2-$dimensional vector space. By cohomology and base change, it even follows that $\mathscr{V}$ is a vector bundle such that $\mathscr{V}_{H}\simeq Ext^{1}(G,H)$, so that the rank of $\mathscr{V}$ is 2.\endproof

We now let$$\mathbb{P}:=\mathbb{P}(\mathscr{V})\stackrel{p}\longrightarrow M_{w}\setminus\{G\},$$which is a $\mathbb{P}^{1}-$bundle over $M_{w}\setminus\{G\}$, and we show that $\mathbb{P}$ is contained in $\widetilde{M}_{v}$:

\begin{prop}
\label{prop:pbundle}Let $(S,v,H)$ be an OLS-triple such that $S$ is an abelian surface. Let $$g:M_{w}\setminus\{G\}\longrightarrow M_{v},\,\,\,\,\,\,\,\,\,\,\,\,g(H):=G\oplus H.$$There is an injective morphism $t:\mathbb{P}\longrightarrow\widetilde{M}_{v}$ such that the following diagram 
$$\begin{array}{ccc}
\mathbb{P} & \stackrel{t}\longrightarrow & \widetilde{M}_{v}\\
\scriptstyle{p}\downarrow & & \downarrow\scriptstyle{\pi_{v}}\\
M_{w}\setminus\{G\} & \stackrel{g}\longrightarrow & M_{v}
\end{array}$$is commutative. 
\end{prop}

\proof We first produce a tautological family $\mathscr{H}$ on $S\times\mathbb{P}$ parameterizing extensions of sheaves $H\in M_{w}\setminus\{G\}$ by $G$. In section 2.2 of \cite{OG2}, and in particular in Proposition 2.2.10, O'Grady gives a modular description of $\widetilde{M}_{v}\setminus\pi_{v}^{-1}(Sing(\Sigma_{v}))$, where $Sing(\Sigma_{v})$ is the singular locus of $\Sigma_{v}$: even if this description is proved for the Mukai vector $(2,0,-2)$, it works in our general setting. Using such a modular description, the tautological family $\mathscr{H}$ induces an injective morphism $t:\mathbb{P}\longrightarrow\widetilde{M}_{v}$ letting the diagram of the statement commute. 

To produce the tautological family $\mathscr{H}$, we use the same procedure described in the Appendix of \cite{Pe}. More precisely, we have the following commutative diagram
$$\begin{array}{ccccc}
S & \stackrel{r_{1}}\longleftarrow & S\times\mathbb{P} & \stackrel{r_{2}}\longrightarrow & \mathbb{P}\\
\scriptstyle{id_{S}}\downarrow & & \scriptstyle{id_{S}\times p}\downarrow & & \downarrow\scriptstyle{p}\\
S & \stackrel{q_{1}}\longleftarrow & S\times (M_{w}\setminus\{G\}) & \stackrel{q_{2}}\longrightarrow & M_{w}\setminus\{G\}
\end{array}$$
and if $\mathscr{T}$ is the tautological bundle of $\mathbb{P}$, then there is a canonical injective morphism $s:\mathscr{T}\longrightarrow p^{*}\mathscr{V}$. Notice that $$p^{*}\mathscr{V}=p^{*}\mathscr{E}xt_{q_{2}}^{1}(q_{1}^{*}G,\mathscr{U}')=\mathscr{E}xt_{r_{2}}^{1}(r_{1}^{*}G,(id_{S}\times p)^{*}\mathscr{U}').$$The morphism $s$ defines then an element $$\sigma\in H^{0}(\mathbb{P},\mathscr{E}xt^{1}_{r_{2}}(r_{1}^{*}G\otimes r_{2}^{*}\mathscr{T},(id_{S}\times p)^{*}\mathscr{U}')).$$Using the spectral sequence associated with the composition of functors, we have and exact sequence $$Ext^{1}(r_{1}^{*}G\otimes r_{2}^{*}\mathscr{T},(id_{S}\times p)^{*}\mathscr{U}')\longrightarrow H^{0}(\mathbb{P},\mathscr{E}xt^{1}_{r_{2}}(r_{1}^{*}G\otimes r_{2}^{*}\mathscr{T},(id_{S}\times p)^{*}\mathscr{U}'))\longrightarrow$$ $$ \longrightarrow H^{2}(\mathbb{P},\mathscr{H}om_{r_{2}}(r_{1}^{*}G\otimes r_{2}^{*}\mathscr{T},(id_{S}\times p)^{*}\mathscr{U}')).$$Now, recall that for every $F\in M_{w}\setminus\{G\}$ we have $Hom(G,F)=0$. It follows from this that $\mathscr{H}om_{r_{2}}(r_{1}^{*}G\otimes r_{2}^{*}\mathscr{T},(id_{S}\times p)^{*}\mathscr{U}')=0$, so that $\sigma$ lifts to an element $\sigma'\in Ext^{1}(r_{1}^{*}G\otimes r_{2}^{*}\mathscr{T},(id_{S}\times p)^{*}\mathscr{U}')$, corresponding to an exact sequence $$0\longrightarrow (id_{S}\times p)^{*}\mathscr{U}'\longrightarrow\mathscr{H}\longrightarrow r_{1}^{*}G\otimes r_{2}^{*}\mathscr{T}\longrightarrow 0.$$The family $\mathscr{H}$ is the one we are looking for, as for every $Q\in\mathbb{P}$, restricting $\mathscr{H}$ to $p^{-1}(Q)$ we get the extension parameterized by $Q$.\endproof

We now state and prove the criterion for divisibility of classes in $H^{2}(\widetilde{M}_{v},\mathbb{Z})$:

\begin{prop}
\label{prop:critdiv}Let $\alpha\in H^{2}(\widetilde{M}_{v},\mathbb{Z})$ and $n\in\mathbb{N}$. There is $\beta\in H^{2}(\widetilde{M}_{v},\mathbb{Z})$ such that $\alpha=n\beta$ if and only if there are $\beta'\in H^{2}(\widetilde{K}_{v},\mathbb{Z})$ and $\beta''\in H^{2}(\mathbb{P},\mathbb{Z})$ such that $\widetilde{j}_{v}^{*}\alpha=n\beta'$ and $t^{*}\alpha=n\beta''$.
\end{prop}

\proof It is immediate to see that if $\alpha=n\beta$ for some $\beta\in H^{2}(\widetilde{M}_{v},\mathbb{Z})$, then $\widetilde{j}_{v}^{*}\alpha=n\widetilde{j}_{v}^{*}\beta$, where $\widetilde{j}_{v}^{*}\beta\in H^{2}(\widetilde{K}_{v},\mathbb{Z})$, and similarly for the restriction to $\mathbb{P}$.

We now show that if $\alpha\in H^{2}(\widetilde{M}_{v},\mathbb{Z})$ is such that there are $\beta'\in H^{2}(\widetilde{K}_{v},\mathbb{Z})$ and $\beta''\in H^{2}(\mathbb{P},\mathbb{Z})$ such that $\widetilde{j}_{v}^{*}\alpha=n\beta'$ and $t^{*}\alpha=n\beta''$, then there is a class $\beta\in H^{2}(\widetilde{M}_{v},\mathbb{Z})$ such that $\alpha=n\beta$. To do so, it is sufficient to show that there is a basis $B=\{\gamma_{1},...,\gamma_{r}\}$ of $H_{2}(\widetilde{M}_{v},\mathbb{Z})/tors$ such that $\alpha\cdot\gamma_{i}\in n\mathbb{Z}$ for every $i=1,...,r$. The following Lemma allows us to construct such a basis $B$.

\begin{lem}
\label{lem:basis}Let $(S,v,H)$ be an OLS-triple such that $S$ is an abelian surface. Let $B_{1}=\{\gamma_{1},...,\gamma_{8}\}$ be a basis of $H_{2}(\widetilde{K}_{v},\mathbb{Z})/tors$ and $B_{2}=\{\delta_{1},...,\delta_{28}\}$ a set of independent elements of $H^{2}(\widetilde{M}_{v},\mathbb{Z})$ such that $\{\widetilde{a}_{v*}\delta_{1},...,\widetilde{a}_{v*}\delta_{28}\}$ is a basis of $H_{2}(S\times\widehat{S},\mathbb{Z})$. Then $B:=\widetilde{j}_{v*}(B_{1})\cup B_{2}$ is a basis of $H_{2}(\widetilde{M}_{v},\mathbb{Z})/tors$. 
\end{lem}

\proof To prove that $B$ is a basis of $H_{2}(\widetilde{M}_{v},\mathbb{Z})/tors$ we construct a basis of $H^{2}(\widetilde{M}_{v},\mathbb{Z})$, and we show that the determinant of the evaluation matrix is 1 or $-1$. Consider the Leray spectral sequence associated to $\widetilde{a}_{v}:\widetilde{M}_{v}\longrightarrow S\times\widehat{S}$, which is
\begin{equation}
\label{eq:leray}
E^{p,q}_{2}=H^{p}(S\times\widehat{S},R^{q}\widetilde{a}_{v*}\mathbb{Z})\Longrightarrow H^{p+q}(\widetilde{M}_{v},\mathbb{Z}).
\end{equation}
Notice that $R^{1}\widetilde{a}_{v*}\mathbb{Z}=0$ as the fibers of $\widetilde{a}_{v}$ are irreducible symplectic manifolds. It follows that $E^{p,1}_{2}=0$ for every $p$, and that $$\widetilde{a}_{v}^{*}:H^{2}(S\times\widehat{S},\mathbb{Z})\longrightarrow H^{2}(\widetilde{M}_{v},\mathbb{Z})$$is injective. Moreover, we have that $E^{0,2}_{2}=H^{0}(S\times\widehat{S},R^{2}\widetilde{a}_{v*}\mathbb{Z})$ is the saturated submodule of $H^{2}(\widetilde{K}_{v},\mathbb{Z})$ given by all the classes which are invariant under monodromy. 

We show that all the classes of $H^{2}(\widetilde{K}_{v},\mathbb{Z})$ are invariant under monodromy, so that $E^{0,2}_{2}=H^{2}(\widetilde{K}_{v},\mathbb{Z})$. Indeed, suppose that $\alpha\in \pi_{v}^{*}(H^{2}(K_{v},\mathbb{Z}))$, and write $\alpha=\pi_{v}^{*}\alpha'$ for some $\alpha'\in H^{2}(K_{v},\mathbb{Z})$: by point 2 of Theorem \ref{thm:pr} we know that $\nu_{v}:v^{\perp}\longrightarrow H^{2}(K_{v},\mathbb{Z})$ is an isomorphism, and as recalled in Remark \ref{oss:lambdatilde} we have a morphism $\lambda_{v}:v^{\perp}\longrightarrow H^{2}(M_{v},\mathbb{Z})$ such that $\nu_{v}=j_{v}^{*}\circ\lambda_{v}$. It follows that the morphism $j_{v}^{*}:H^{2}(M_{v},\mathbb{Z})\longrightarrow H^{2}(K_{v},\mathbb{Z})$ is surjective, hence there is $\alpha''\in H^{2}(M_{v},\mathbb{Z})$ such that $\alpha'=j_{v}^{*}(\alpha'')$. Then $\alpha=\widetilde{j}_{v}^{*}(\pi_{v}^{*}(\alpha''))$, and $\pi_{v}^{*}(H^{2}(K_{v},\mathbb{Z}))\subseteq E^{0,2}_{2}$.

Now, by point 2 of Theorem \ref{thm:w2hs} we know that $H^{2}(\widetilde{K}_{v},\mathbb{Z})$ is generated by the pull-back of $H^{2}(K_{v},\mathbb{Z})$ and by $c_{1}(A_{v})$, where $A_{v}\in Pic(\widetilde{K}_{v})$ is such that $A_{v}^{\otimes 2}=\mathscr{O}(\widetilde{\Sigma}_{v})$. Now, notice that $c_{1}(\widetilde{\Sigma}_{v})=\widetilde{j}_{v}^{*}(c_{1}(\widetilde{E}))$, so that $c_{1}(\widetilde{\Sigma}_{v})\in E^{0,2}_{2}$. As $E^{0,2}_{2}$ is a saturated submodule of $H^{2}(\widetilde{K}_{v},\mathbb{Z})$, it follows that $c_{1}(A_{v})\in E^{0,2}_{2}$. In conclusion, we have $E^{0,2}_{2}=H^{2}(\widetilde{K}_{v},\mathbb{Z})$.

Consider now the Leray spectral sequence $$H^{p}(S\times\widehat{S},R^{q}\widetilde{a}_{v*}\mathbb{Q})\Longrightarrow H^{p+q}(\widetilde{M}_{v},\mathbb{Q}),$$which degenerates by Deligne's Theorem. We have then an exact sequence $$0\longrightarrow H^{2}(S\times\widehat{S},\mathbb{Q})\longrightarrow H^{2}(\widetilde{M}_{v},\mathbb{Q})\longrightarrow H^{2}(\widetilde{K}_{v},\mathbb{Q})\longrightarrow 0,$$so that $b_{2}(\widetilde{M}_{v})=b_{2}(S\times\widehat{S})+b_{2}(\widetilde{K}_{v})=36$. Now, the Leray spectral sequence (\ref{eq:leray}) gives us an exact sequence $$0\longrightarrow H^{2}(S\times\widehat{S},\mathbb{Z})\stackrel{\widetilde{a}_{v}^{*}}\longrightarrow H^{2}(\widetilde{M}_{v},\mathbb{Z})\stackrel{\widetilde{j}_{v}^{*}}\longrightarrow H^{2}(\widetilde{K}_{v},\mathbb{Z})\longrightarrow H^{3}(S\times\widehat{S},\mathbb{Z}).$$As $b_{2}(\widetilde{M}_{v})=b_{2}(S\times\widehat{S})+b_{2}(\widetilde{K}_{v})$, it follows that the image of the morphism $H^{2}(\widetilde{K}_{v},\mathbb{Z})\longrightarrow H^{3}(S\times\widehat{S},\mathbb{Z})$ is contained in the torsion part of $H^{3}(S\times\widehat{S},\mathbb{Z})$, which is surely trivial. 

We finally get an exact sequence $$0\longrightarrow H^{2}(S\times\widehat{S},\mathbb{Z})\stackrel{\widetilde{a}_{v}^{*}}\longrightarrow H^{2}(\widetilde{M}_{v},\mathbb{Z})\stackrel{\widetilde{j}_{v}^{*}}\longrightarrow H^{2}(\widetilde{K}_{v},\mathbb{Z})\longrightarrow 0.$$It follows that there is a set $\{\zeta_{1},...,\zeta_{8}\}$ of elements of $H^{2}(\widetilde{M}_{v},\mathbb{Z})$ which restricts to a basis of $H^{2}(\widetilde{K}_{v},\mathbb{Z})$, and if $\{\xi_{1},...,\xi_{28}\}$ is a basis of $H^{2}(S\times\widehat{S},\mathbb{Z})$ then the set $B':=\{\widetilde{a}_{v}^{*}\xi_{1},...,\widetilde{a}_{v}^{*}\xi_{28},\zeta_{1},...,\zeta_{8}\}$ is a basis of $H^{2}(\widetilde{M}_{v},\mathbb{Z})$.

To show that $B$ is a basis for $H_{2}(\widetilde{M}_{v},\mathbb{Z})/tors$, we show that the evaluation matrix $M$ (whose entries are the evaluations of the elements of $B'$ on those of $B$) has determinant 1 or $-1$. Notice that $M$ has four blocks, which are $M_{1}:=[\zeta_{i}(\gamma_{j})]$, $M_{2}=[\zeta_{i}(\delta_{j})]$, $M_{3}=[\widetilde{a}_{v}^{*}\xi_{i}(\gamma_{j})]$ and $M_{4}=[\widetilde{a}_{v}^{*}\xi_{i}(\delta_{j})]$. By definition, we have $|det(M_{1})|=|det(M_{4})|=1$, and $M_{3}=0$, so that in conclusion $|det(M)|=1$, and we are done.\endproof

We now conclude the proof of Proposition \ref{prop:critdiv}. We construct a basis $B$ of $H_{2}(\widetilde{M}_{v},\mathbb{Z})/tors$ as in the statement of Lemma \ref{lem:basis}, and such that $\alpha(B)\subseteq n\mathbb{Z}$. Let $\{\gamma_{1},...,\gamma_{8}\}$ be a basis of $H_{2}(\widetilde{K}_{v},\mathbb{Z})/tors$: as $\widetilde{j}_{v}^{*}\alpha=n\beta'$ for $\beta'\in H^{2}(\widetilde{K}_{v},\mathbb{Z})$ it follows that $$\alpha(\widetilde{j}_{v*}\gamma_{i})=\widetilde{j}_{v}^{*}\alpha(\gamma_{i})=n\beta'(\gamma_{i})\in n\mathbb{Z}$$for $i=1,...,8$. It remains to find $\{\delta_{1},...,\delta_{28}\}\subseteq H^{2}(\widetilde{M}_{v},\mathbb{Z})$ which surjects onto a basis of $H^{2}(S\times\widehat{S},\mathbb{Z})$, and such that $\alpha(\delta_{i})\in n\mathbb{Z}$ for $i=1,...,28$.

To do so, notice that we have a commutative diagram
\begin{equation}
\label{eq:cmp}
\begin{array}{ccc}
\mathbb{P} & \stackrel{t}\longrightarrow & \widetilde{M}_{v}\\
\scriptstyle{p}\downarrow & & \downarrow\scriptstyle{\widetilde{a}_{v}}\\
M_{w}\setminus\{G\} & \stackrel{a_{w}}\longrightarrow & S\times\widehat{S}
\end{array}
\end{equation}
where by abuse of notation we write $a_{w}$ for the restriction to $M_{w}\setminus\{G\}$ of the isomorphism $a_{w}:M_{w}\longrightarrow S\times\widehat{S}$. Since $\mathbb{P}$ is the projectivized of a rank 2 vector bundle $\mathscr{V}$, from the Leray spectral sequence of the fibration $p$ the morphism $p^{*}:H^{2}(M_{w}\setminus\{G\},\mathbb{Z})\longrightarrow H^{2}(\mathbb{P},\mathbb{Z})$ is a saturated inclusion. It follows that $$p_{*}:H_{2}(\mathbb{P},\mathbb{Z})\longrightarrow H_{2}(M_{w}\setminus\{G\},\mathbb{Z})$$is surjective. Let $\{\delta'_{1},...,\delta'_{28}\}$ be a set of elements of $H_{2}(\mathbb{P},\mathbb{Z})$ which surjects onto a basis of $H_{2}(M_{w}\setminus\{G\},\mathbb{Z})$, and let $\delta_{i}:=t_{*}\delta'_{i}$ for $i=1,...,28$. As $a_{w}$ induces an isomorphism between $H^{2}(M_{w}\setminus\{G\},\mathbb{Z})$ and $H^{2}(S\times\widehat{S},\mathbb{Z})$, by the commutativity of the diagram (\ref{eq:cmp}) we then see that $\{\delta_{1},...,\delta_{28}\}$ is a set of elements of $H_{2}(\widetilde{M}_{v},\mathbb{Z})$ which surjects onto a basis of $H_{2}(S\times\widehat{S},\mathbb{Z})$. Moreover, as $t^{*}\alpha=n\beta''$ for some $\beta''\in H^{2}(\mathbb{P},\mathbb{Z})$, we then have $$\alpha(\delta_{i})=t^{*}\alpha(\delta'_{i})=n\beta''(\delta'_{i})\in n\mathbb{Z},$$and we are done.\endproof

\begin{oss}
\label{oss:36}{\rm In the proof of Lemma \ref{lem:basis} we showed that $b_{2}(\widetilde{M}_{v})=36$, and that $H^{2}(\widetilde{M}_{v},\mathbb{Z})$ is a free $\mathbb{Z}-$module.}
\end{oss}

\paragraph{Conclusion of the proof of Theorem \ref{thm:mainab}.}

We are now ready to prove the remaining part of Theorem \ref{thm:mainab}, which is the following:

\begin{prop}
\label{prop:mvnlf}Let $(S,v,H)$ be an OLS-triple such that $S$ is an abelian surface. Then $M_{v}(S,H)$ is $2-$factorial.
\end{prop}

\proof In order to prove the statement, by Proposition \ref{prop:mvlf2f} we just need to show that $Pic(M_{v})$ is strictly contained in $Pic(M_{v}^{s})$. This is the case if there is a line bundle on $\widetilde{M}_{v}$ whose restriction to $\widetilde{K}_{v}$ is $A_{v}$. We show that it exists by presenting a line bundle $L\in Pic(\widetilde{M}_{v})$ whose restriction to $\widetilde{K}_{v}$ is $\mathscr{O}(\widetilde{\Sigma}_{v})$, and such that there is $\beta\in H^{2}(\widetilde{M}_{v},\mathbb{Z})$ such that $c_{1}(L)=2\beta$. 

First of all, recall that we defined a rank 2 vector bundle $\mathscr{V}$ on $M_{w}\setminus\{G\}$ such that $\mathbb{P}=\mathbb{P}(\mathscr{V})$. As $a_{w}:M_{w}\setminus\{G\}\longrightarrow S\times\widehat{S}$ is an open embedding, there is a line bundle $\mathscr{L}\in Pic(S\times\widehat{S})$ such that $a_{w}^{*}(\mathscr{L})=det(\mathscr{V})$.

We define $L:=\widetilde{E}+\widetilde{a}_{v}^{*}(\mathscr{L})$: notice that $$\widetilde{j}_{v}^{*}(L)=\widetilde{j}_{v}^{*}(\widetilde{E})=\mathscr{O}(\widetilde{\Sigma}_{v}).$$By Proposition \ref{prop:critdiv}, it is then sufficient to show that $\widetilde{j}_{v}^{*}(L)$ and $t^{*}(L)$ admit square roots in $Pic(\widetilde{K}_{v})$ and $Pic(\mathbb{P})$ respectively. We already know that $\widetilde{j}_{v}^{*}(L)=A_{v}^{\otimes 2}$, hence it remains to show that $t^{*}(L)$ admits a square root in $Pic(\mathbb{P})$.

Now, by adjunction $$t^{*}(\widetilde{E})=K_{\mathbb{P}/M_{w}\setminus\{G\}},$$the relative canonical divisor of the $\mathbb{P}^{1}-$fibration $p:\mathbb{P}\longrightarrow M_{w}\setminus\{G\}$. By the Euler exact sequence, as $\mathbb{P}=\mathbb{P}(\mathscr{V})$ and $\mathscr{V}$ is a rank 2 vector bundle over $M_{w}\setminus\{G\}$, we have that $$K_{\mathbb{P}/M_{w}\setminus\{G\}}=\mathscr{O}_{\mathbb{P}}(2)-p^{*}det(\mathscr{V}).$$Using the commutativity of the diagram (\ref{eq:cmp}), this implies that $$t^{*}(L)=t^{*}(\widetilde{E})+t^{*}(\widetilde{a}_{v}^{*}(\mathscr{L}))=K_{\mathbb{P}/M_{w}\setminus\{G\}}+p^{*}det(\mathscr{V})=\mathscr{O}_{\mathbb{P}}(2).$$Hence $\mathscr{O}_{\mathbb{P}}(1)$ is a square root of $t^{*}(L)$ in $Pic(\mathbb{P})$.\endproof

\section{Appendix}

In this paper we studied the index of factoriality of $K_{v}(S,H)$ for an OLS-triple $(S,v,H)$ where $S$ is abelian. It seems natural to us to provide a complete description of the index of factoriality of $K_{v}(S,H)$: this is the main motivation for this Appendix. Before starting, we recall some general facts about $K_{v}(S,H)$: first, $K_{v}(S,H)$ is the fiber over $(0,\mathscr{O}_{S})$ of the morphism $a_{v}:M_{v}(S,H)\longrightarrow S\times\widehat{S}$ we defined in the introduction. If $v^{2}>0$, then this map is dominant, it is an isotrivial fibration, and there is a finite étale morphism $$\tau:K_{v}\times S\times\widehat{S}\longrightarrow M_{v},\,\,\,\,\,\,\,\,\,\,\,\,\tau(F,p,L):=t_{p}^{*}(F)\otimes L,$$where $t_{p}:S\longrightarrow S$ is the translation via $p$.

\begin{oss}
\label{oss:irrkv}
{\rm If $v^{2}>0$, then $K_{v}(S,H)$ is normal and irreducible. The normality of $K_{v}$ follows from the normality of $M_{v}$ (see Theorem 4.4 of \cite{KLS}) and from the finite étale cover $\tau$. The irreducibility is a consequence of the irreducibility of $M_{v}$ and of the one of $K_{w}$. Indeed, as $K_{v}$ is normal, to show that it is irreducible it is sufficient to show that it is connected. 

Now, let $E\in K_{v}$: there are $E'_{0}\in M_{w}$ and an irreducible curve $D'\subseteq M_{v}$ passing through $(E'_{0})^{\oplus m}$ and $E$, as $M_{v}$ is irreducible (by Theorem 4.4 of \cite{KLS}). Let $D''$ be an irreducible component of the curve $\tau^{-1}(D')$ passing through $(E,0,\mathscr{O}_{S})$. Moreover, let $p_{1}:K_{v}\times S\times\widehat{S}\longrightarrow K_{v}$ be the projection, and $D:=p_{1}(D'')$. This is an irreducible curve in $K_{v}$ passing through $E$ and $E_{0}^{\oplus m}$, where $E_{0}:=t_{p}^{*}F_{0}\otimes L$ for some $F_{0}\in M_{w}$, $p\in S$ and $L\in\widehat{S}$. 

In conclusion, in order to show that $K_{v}$ is connected, it is sufficient to show that the locus $Y$ of $K_{v}$ parameterizing direct sums of stable sheaves of Mukai vector $w$ is connected. Notice that we have a surjective map from the locus $Y'\subseteq M_{w}^{m}$ of $m-$tuples $(F_{1},...,F_{m})\in M_{w}^{m}$ such that $\sum_{i=1}^{m}a_{w}(F_{i})=(0,\mathscr{O}_{S})$ to $Y$. It is then sufficient to show that $Y'$ is connected. 

Now, notice that $Y'$ is a projective variety. Moreover, let $s:(S\times\widehat{S})^{m}\longrightarrow S\times\widehat{S}$ be the sum morphism: we have a fibration $f:Y'\longrightarrow s^{-1}(0,\mathscr{O}_{S})$ with fibers isomorphic to $K_{w}^{m}$. As $s^{-1}(0,\mathscr{O}_{S})$ is irreducible, and $K_{w}$ are irreducible by \cite{Y5}, it follows that $Y'$ is connected.}
\end{oss}

The only results we have up to now about the index of factoriality of $K_{v}(S,H)$ are the following:
\begin{enumerate}
 \item if $m=1$, then $K_{v}(S,H)$ is smooth, and hence locally factorial;
 \item if $m=2$ and $w^{2}=2$, then $K_{v}(S,H)$ is $2-$factorial by Theorem \ref{thm:mainab};
\end{enumerate}
In order to complete the description of the index of factoriality of $K_{v}$, we need to study the two following cases:
\begin{enumerate}
 \item $m\geq 2$ and $w^{2}=0$;
 \item $m=2$ and $w^{2}\geq 4$, or $m\geq 3$ and $w^{2}\geq 2$.
\end{enumerate} 

We begin with the index of factoriality when $m\geq 2$ and $w^{2}=0$. In this case, the moduli space $M_{v}(S,H)$ is isomorphic to $M_{w}^{(m)}$, the $m-$th symmetric product of the abelian surface $M_{w}$. We have the sum morphism $s:M_{w}^{(m)}\longrightarrow M_{w}$, and $K_{v}\simeq s^{-1}(0)$. Then $K_{v}$ is $2-$factorial, as the following result shows:

\begin{prop}
\label{prop:symmab}
Let $S$ be an abelian surface, $n\in\mathbb{N}$. Let $s:S^{(n+1)}\longrightarrow S$ be the sum morphism, and $K^{n}_{sing}(S):=s^{-1}(0)$. The variety $K^{n}_{sing}(S)$ is $2-$factorial.
\end{prop}

\proof Let $\rho_{n+1}:Hilb^{n+1}(S)\longrightarrow S^{(n+1)}$ be the Hilbert-Chow morphism, and $K^{n}(S):=\rho_{n+1}^{-1}(K^{n}_{sing}(S))$: if $n=1$, it is the Kummer surface associated to $S$; if $n\geq 2$, it is a generalized Kummer variety of dimension $2n$.

Let us first show that $K^{n}_{sing}(S)$ is $2-$factorial if $n\geq 2$. By abuse of notation, we still write $\rho_{n+1}:K^{n}(S)\longrightarrow K^{n}_{sing}(S)$. Let $E$ be the exceptional divisor of $\rho_{n+1}$ on $K^{n}(S)$, and $A\in Pic(K^{n}(S))$ be such that $A^{\otimes 2}=\mathscr{O}(E)$.

Lemma 2.1 of \cite{Y3} and Proposition 8 of \cite{Bea} show that we have an equality $H^{2}(K^{n}(S),\mathbb{Z})= j(H^{2}(S,\mathbb{Z}))\oplus\mathbb{Z}\cdot A$, where $j:H^{2}(S,\mathbb{Z})\longrightarrow H^{2}(K^{n}(S),\mathbb{Z})$ is an injective morphism. It follows that $$Pic(K^{n}(S))= j(NS(S))\oplus\mathbb{Z}\cdot A.$$The morphism $j$ works as follows (see the proof of Proposition 8 of \cite{Bea}): let $\beta\in NS(S)$, and let $p_{i}:S^{n+1}\longrightarrow S$ be the projection to the $i-$th factor. Consider the class $\beta':=\sum_{i=1}^{n+1}p_{i}^{*}\beta\in H^{2}(S^{n+1},\mathbb{Z})$. By Lemma 6.1 of \cite{F} there is $L\in Pic(S^{(n+1)})$ such that $\beta'=q^{*}c_{1}(L)$, where $q:S^{n+1}\longrightarrow S^{(n+1)}$ is the quotient morphism. Let now $k:K^{n}(S)\longrightarrow Hilb^{n+1}(S)$ be the inclusion morphism: then $j(\beta):=k^{*}\rho_{n+1}^{*}(L)$. 

If $h:K^{n}_{sing}(S)\longrightarrow S^{(n+1)}$ is the inclusion, we then have $j(\beta)=\rho_{n+1}^{*}h^{*}(L)$, so that $j(NS(S))=\rho_{n+1}^{*}(Pic(K^{n}_{sing}(S)))$. By the exact sequence $$0\longrightarrow\mathbb{Z}\stackrel{\iota}\longrightarrow Pic(K^{n}(S))\stackrel{r}\longrightarrow Pic(K^{n}(S)\setminus E)\longrightarrow 0,$$where $\iota(1):=\mathscr{O}(E)$ and $r$ is the restriction morphism, we get $$Pic(K^{n}(S)\setminus E)=Pic(K^{n}(S))/\mathbb{Z}\cdot\mathscr{O}(E)=j(NS(S))\oplus\mathbb{Z}/2\mathbb{Z}\cdot A.$$By Remark \ref{oss:cw} we have $A^{1}(K^{n}_{sing}(S))=Pic(K^{n}(S)\setminus E)$, so that $$A^{1}(K^{n}_{sing}(S))/Pic(K^{n}_{sing}(S))= (j(NS(S))\oplus\mathbb{Z}/2\mathbb{Z}\cdot A)/j(NS(S))\simeq\mathbb{Z}/2\mathbb{Z},$$and the $2-$factoriality of $K^{n}_{sing}(S)$ is shown. 

There is still the case of $K^{1}(S)$, which is the Kummer surface associated to $S$, and $K^{1}_{sing}(S)$ is the singular Kummer surface. The $2-$factoriality of $K^{1}_{sing}(S)$ can be proved either by using an argument similar to the one we used here (where one replaces the irreducible divisor $E$ by the 16 nodal curves), or by recalling that $K^{1}_{sing}(S)$ is the quotient of $S$ by the action of $-1$: we leave the details to the reader.\endproof

We are left with the determination of the index of factoriality of $K_{v}(S,H)$ for $v=mw$, where either $m=2$ and $w^{2}\geq 4$, or $m\geq 3$ and $w^{2}\geq 2$. We show that $K_{v}(S,H)$ is locally factorial in this case, by strictly following the argument used by Kaledin, Lehn and Sorger in \cite{KLS} to prove the local factoriality of $M_{v}(S,H)$ under the same hypothesis. Before proving the Proposition, we give the following definition:

\begin{defn}
Let $F$ and $F'$ be two $H-$polystable sheaves on $S$ with Mukai vector $v=mw$. Write $F=\bigoplus_{i=1}^{s}F_{i}^{\oplus m_{i}}$ and $F'=\bigoplus_{i=1}^{s'}(F'_{i})^{\oplus m'_{i}}$, where $F_{i}$ (resp. $F'_{j}$) is $H-$stable and $F_{i}\nsimeq F_{k}$ (resp. $F'_{i}\nsimeq F'_{k}$) for every $i,k=1,...,s$ (resp. $j,k=1,...,s'$). We say that $F$ and $F'$ \textit{have the same polystable structure} if $s=s'$, $v(F_{i})=v(F'_{i})$ and $m_{i}=m'_{i}$ for every $i=1,...,s$.
\end{defn}

We now prove the local factoriality of $K_{v}$:

\begin{prop}
\label{prop:kvlf}
Let $S$ be an abelian surface, $v=mw$ a Mukai vector and $H$ a $v-$generic polarization on $S$. Suppose either that $m=2$ and $w^{2}\geq 4$, or that $m\geq 3$ and $w^{2}\geq 2$. Then $K_{v}(S,H)$ is locally factorial.
\end{prop}

\proof We recall some elements of the construction of $M_{v}$. First, there are $k,N\in\mathbb{N}$ such that every $H-$semistable sheaf of Mukai vector $v$ is quotient of $\mathscr{H}:=\mathscr{O}_{S}(-kH)^{\oplus N}$. Let $Q_{v}$ be the Quot scheme parameterizing quotients of $\mathscr{H}$ with Mukai vector $v$, and let $R_{v}\subseteq Q_{v}$ be the open subset parameterizing $H-$semistable quotients. The group $PGL(N)$ acts on $R_{v}$, and the quotient is $M_{v}$: let $p:R_{v}\longrightarrow M_{v}$ be the quotient morphism, and let $R_{v}^{0}:=p^{-1}(K_{v})$.

As a first step, following \cite{KLS} we show that $R_{v}^{0}$ is locally factorial by showing that it is locally complete intersection (we will write l.c.i. in the following) and regular in codimension 3, i. e. $R_{v}^{0}$ has the property ($R_{3}$). 

By \cite{KLS} we know that $R_{v}$ is l.c.i. As $R_{v}^{0}=p^{-1}(a_{v}^{-1}(0,\mathscr{O}_{S}))$, and $(0,\mathscr{O}_{S})$ is a smooth point in $S\times\widehat{S}$, it follows that $R_{v}^{0}$ is l.c.i.

We now show that $R_{v}^{0}$ has the property ($R_{3}$) by showing that its singular locus has codimension at least 4. Now, let $F$ be an $H-$semistable sheaf with Mukai vector $v$ on $S$, and consider a family $\mathscr{F}$ on $S\times S\times\widehat{S}$ such that for every $(p,L)\in S\times\widehat{S}$ we have $\mathscr{F}_{|p_{23}^{-1}(p,L)}=t_{p}^{*}(F)\otimes L$, as constructed in section 4.2.2. This family induces a modular morphism $$f:S\times\widehat{S}\longrightarrow M_{v},$$such that $a_{v}\circ f:S\times\widehat{S}\longrightarrow S\times\widehat{S}$ is étale by section 1.2 of \cite{Y3}. The sheaf $\mathscr{H}om_{p_{23}}(p_{1}^{*}\mathscr{H},\mathscr{F})$ is a vector bundle on $S\times\widehat{S}$, and let $$\pi:\mathbb{P}(\mathscr{H}om_{p_{23}}(p_{1}^{*}\mathscr{H},\mathscr{F}))\longrightarrow S\times\widehat{S}$$be the associated  projective bundle. Let $U\subseteq\mathbb{P}(\mathscr{H}om_{p_{23}}(p_{1}^{*}\mathscr{H},\mathscr{F}))$ be the open subset parameterizing surjective morphisms: there is a natural injective morphism $g:U\longrightarrow R_{v}$ such that the following diagram is commutative
$$\begin{array}{ccc}
U & \stackrel{g}\longrightarrow & R_{v}\\
\scriptstyle{\pi}\downarrow & & \downarrow\scriptstyle{p}\\
S\times\widehat{S} & \stackrel{f}\longrightarrow & M_{v}
\end{array}$$
Let $[q:\mathscr{H}\longrightarrow E]\in U$ be a quotient of $\mathscr{H}$. The morphism $$d_{[q]}\pi:T_{[q]}U\longrightarrow T_{\pi([q])}(S\times\widehat{S})$$is surjective: it follows that $$d_{[q]}(a_{w}\circ p):T_{[q]}R_{v}\longrightarrow T_{a_{w}(E)}(S\times\widehat{S})$$is surjective. This implies that if $[q]$ is a smooth point of $R_{v}$ which lies in $R_{v}^{0}$, then $[q]$ is a smooth point of $R_{v}^{0}$. As a consequence, the singular locus $Sing(R_{v}^{0})$ of $R_{v}^{0}$ is contained in the singular locus $Sing(R_{v})$ of $R_{v}$.

Now, $R_{v}$ is regular in codimension 3 by \cite{KLS}. Moreover, if $[F]\in M_{v}$ is a polystable sheaf, the dimension of $Sing(R_{v})\cap p^{-1}([F])$ depends only on the polystable structure of $F$. Similarily, if $(p,L)\in S\times\widehat{S}$, then the dimension of the locus of the polystable sheaves in $a_{v}^{-1}(p,L)$ having the same polystable structure does not depend on $(p,L)$. As $Sing(R_{v}^{0})\subseteq Sing(R_{v})$, it follows that $R_{v}^{0}$ is regular in codimension 3.

As $R_{v}^{0}$ is locally complete intersection and regular in codimension 3, by Corollaire 3.14 of Exposé XI in \cite{G} we have that $R_{v}^{0}$ is locally factorial. As remarked in \cite{KLS}, this allows us to apply Théorème A of \cite{D2} to show that $K_{v}$ is locally factorial: $K_{v}$ is the quotient of $R_{v}^{0}$ by the action of $PGL(N)$, hence $K_{v}$ is locally factorial at a point $E\in K_{v}\setminus K_{v}^{s}$ if and only if the isotropy subgroup $PAut(E)$ of any point $[q]$ in the closed orbit $p^{-1}(E)\subseteq R_{v}^{0}$ acts trivially on the fiber $L_{[q]}$ for every $PGL(N)-$linearized line bundle $L$ on an invariant open neighbourhood of the orbit of $[q]$.

First, consider a point $E_{0}^{\oplus m}\in K_{v}$ for some $E_{0}\in M_{w}$. The isotropy subgroup of any point $[q]\in p^{-1}(E_{0}^{\oplus m})$ is $PGL(m)$, and hence has no non-trivial characters (see Corollary 5.2 in \cite{KLS}). The point $E_{0}^{\oplus m}$ is then locally factorial.

Let now $E=\bigoplus_{i=1}^{s}E_{i}^{\oplus n_{i}}\in K_{v}$, where $E_{i}$ is stable, $v(E_{i})=m_{i}w$ for some $m_{i}\in\mathbb{N}$, and $m=\sum_{i=1}^{s}m_{i}n_{i}$. Let $[q]\in p^{-1}(E)$, and as $E\in K_{v}$ we have $[q]\in R_{v}^{0}$. Consider the moduli space $M_{m_{i}w}$ for $i=1,...,s$: as every sheaf $F\in M_{m_{i}w}$ is $H-$semistable, there is an integer $N_{i}\in\mathbb{N}$ such that $F$ is quotient of $\mathscr{H}_{i}:=\mathscr{O}_{S}(-kH)^{\oplus N_{i}}$. We suppose that every sheaf parameterized by $M_{v}$ is quotient of $\mathscr{H}=\bigoplus_{i=1}^{s}\mathscr{H}_{i}^{\oplus n_{i}}$. We let $Q_{m_{i}w}$ be the Quot scheme parameterizing quotients of $\mathscr{H}_{i}$, and $R_{m_{i}w}$ the open subset of $Q_{m_{i}w}$ parameterizing $H-$semistable quotients. Now, consider the morphism $$\Phi:\prod_{i=1}^{s}R_{m_{i}w}\longrightarrow R_{v},\,\,\,\,\,\,\,\,\Phi([q_{1}],...,[q_{s}]):=\bigg[\bigoplus_{i=1}^{s}q_{i}^{\oplus n_{i}}\bigg].$$As shown in \cite{KLS}, $M_{m_{i}w}$ and $R_{m_{i}w}$ are irreducible. It follows that the image $Z$ of the morphism $\Phi$ is irreducible. By definition, we have that $[q]\in Z$.

Now, consider the group $G:=(\prod_{i=1}^{s}GL(n_{i}))/\mathbb{C}^{*}\subseteq PGL(N)$. This group fixes $Z\cap R_{v}^{0}$ pointwise, and it is the stabilizer of $[q]$. Moreover, if $L$ is any $PGL(N)-$linearized line bundle on $R^{0}_{v}$, then $G$ acts on $L_{Z\cap R^{0}_{v}}$ via a locally constant character. 

We claim that there is a connected curve $C\subseteq Z\cap R^{0}_{v}$ passing through $[q]$ and through a point $[q']\in p^{-1}(E_{0}^{\oplus m})$ for some $E_{0}^{\oplus m}\in K_{v}$. Once this claim is proved, we conclude the proof: as $C$ is connected, the action of $G$ on $L_{|C}$ is via a constant character. As $G$ is contained in the stabilizer of $[q']$, its action on $L_{[q']}$ has to be trivial: this implies that the action of $G$ on $L_{[q]}$ is trivial, and we are done by Théorème A of \cite{D2}.

We are left with the proof of the claim. Consider the morphism $$h:\prod_{i=1}^{s}M_{m_{i}w}\longrightarrow M_{v},\,\,\,\,\,\,\,\,\,\,h(F_{1},...,F_{s}):=\bigoplus_{i=1}^{s}F_{i}^{\oplus n_{i}},$$and the commutative diagram 
$$\begin{array}{ccc}
\prod_{i=1}^{s}R_{m_{i}w} & \stackrel{\Phi}\longrightarrow & R_{v}\\
\scriptstyle{g}\downarrow & & \downarrow\scriptstyle{p}\\
\prod_{i=1}^{s}M_{m_{i}w} & \stackrel{h}\longrightarrow & M_{v}
\end{array}$$
Let $Q=([q_{1}],...,[q_{s}])\in\prod_{i=1}^{s}R_{m_{i}w}$ be such that $\Phi(Q)=[q]$.

First, following the same construction as in Remark \ref{oss:irrkv}, there is an irreducible curve $D\subseteq K_{v}$ passing through $E$ and a point $E_{0}^{\oplus m}\in K_{v}$ for some $E_{0}\in M_{w}$: this is constructed as $p_{1}(\tau^{-1}(D''))$, where $p_{1}:K_{v}\times S\times\widehat{S}\longrightarrow K_{v}$ is the projection, and $D''$ is an irreducible component of $\tau^{-1}(D')$ for some irreducible curve passing through $E$ and a point $(E'_{0})^{m}\in M_{v}$. Notice that as $\tau$ preserves the polystable structure, then $D$ is contained in the image of $h$.

Now, let $[q']\in p^{-1}(E_{0}^{\oplus m})$. As $D$ is contained in the image of $h$, and $h$ is étale on its image on a neighbourhood of $g(Q)$, there is an irreducible component $D_{1}$ of $h^{-1}(D)$ passing through $g(Q)$ and a point of $h^{-1}(E_{0}^{\oplus m})$. Moreover, $D_{1}$ is the quotient of $g^{-1}(D_{1})$ by a connected group: as $D_{1}$ is connected, it follows that $g^{-1}(D_{1})$ is connected. There is then a connected curve $D_{2}\subseteq g^{-1}(D_{1})$ passing through $Q$ and a point of $\Phi^{-1}([q'])$. We finally define $C:=\Phi(D_{2})$: this is a connected curve contained in $Z\cap R_{v}^{0}$, which passes through $[q]$ and $[q']$. This concludes the proof of the claim.\endproof

\subsection*{Acknowledgements}We thank Christoph Sorger and Olivier Serman for useful discussions. The first author was supported by the SFB/TR 45 'Periods, Moduli spaces and Arithmetic of Algebraic Varieties' of the DFG (German Research Fundation).

\par\bigskip
\par\bigskip
Arvid Perego, Institut Elie Cartan de Nancy (UMR 7502 du CNRS), Université Henri Poincaré Nancy 1, B.P. 70239, 54506 Vandoeuvre-lès-Nancy Cedex, France.

\textit{E-mail address:} \texttt{Arvid.Perego@iecn.u-nancy.fr}
\par\bigskip
Antonio Rapagnetta, Dipartimento di Matematica dell'Universit\`a di Roma II - Tor Vergata, 00133 Roma, Italy.

\textit{E-mail address:} \texttt{rapagnet@mat.uniroma2.it} 


\begin{thebibliography}{10}
\addcontentsline{toc}{chapter}{Bibliography}
\bibitem{Bea} Beauville, A.: Variétés de K\"ahler dont la première classe de Chern est nulle. J. Differential Geom. Volume 18, Number 4 (1983), 755-782.
\bibitem{BLS} Beauville, A., Laszlo, Y., Sorger, C.: The Picard group of the moduli of $G-$bundles on a curve. Compositio Math. \textbf{112}, 183--216 (1998).
\bibitem{BGS} Boissière, S., Gabber, O., Serman, O.: Sur le produit de variétés localement factorielles ou $\mathbb{Q}-$factorielles. Preprint version: arXiv:1104.1861.
\bibitem{BK} Boysal, A., Kumar, S.: Explicit Determination of the Picard Group of Moduli Spaces of Semi-Stable $G-$Bundles on Curves. Math. Ann. \textbf{332}, 4, 823--842, (2005).
\bibitem{D1} Dr\'ezet, J.-M.: Groupe de Picard des vari\'et\'es de modules de faisceaux semi-stables sur $\mathbb{P}_{2}(\mathbb{C})$. Ann. Ist. Fourier \textbf{38}, 105--168 (1988).
\bibitem{D2} Dr\'ezet, J.-M.: Points non factoriels des vari\'et\'es de modules de faisceaux semi-stables sur une surface rationnelle. J. reine und angew. Math. \textbf{290}, 99--127 (1991).
\bibitem{DN} Dr\'ezet, J.-M., Narasimhan, M. S.: Groupe de Picard des vari\'et\'es de modules de fibres semistables sur les courbes algebriques, Invent. Math. \textbf{97}, 1, 53--94 (1989).
\bibitem{F} Fogarty, J.: Algebraic Families on an Algebraic Surface, II, the Picard Scheme of the Punctual Hilbert Scheme, American Journal of Mathematics \textbf{95}, 3, 660--687 (1973).
\bibitem{G} Grothendieck, A.: Cohomologie locale des faisceaux cohérents et Théorèmes de Lefschetz locaux et globaux. Séminaire de Géométrie Algébrique du Bois-Marie 1962 (SGA 2). North Holland Publishing Company, Amsterdam (1968). 
\bibitem{H} Hartshorne, R.: Algebraic Geometry. GTM \textbf{52}, Springer-Verlag, New York (1977).
\bibitem{KLS} Kaledin, D., Lehn, M., Sorger, C.: Singular Symplectic Moduli Spaces. Invent. Math. \textbf{164}, 591--614 (2006).
\bibitem{La} Laszlo, Y.: About $G-$bundles over elliptic curves. Ann. Inst. Fourier \textbf{48}, 2, 413--424 (1998). 
\bibitem{LaS} Laszlo, Y., Sorger, C.: The line bundles on the moduli of parabolic $G-$bundles over curves and their sections.
Ann. Sci. \'Ecole Norm. Sup. \textbf{30}, 4, 499--525 (1997)
\bibitem{LS} Lehn, M., Sorger, C.: La singularit\'e de O'Grady. J. Algebraic Geometry \textbf{15}, 756--770 (2006).
\bibitem{M1} Mukai, S.: Symplectic structure of the moduli space of sheaves on an abelian or K3 surface, Invent. Math. \textbf{77}, 101--116 (1984).
\bibitem{M2} Mukai, S.: Duality between $D(X)$ and $D(\widehat{X})$ with its application to Picard sheaves, Nagoya Math. J. \textbf{81}, 153--175 (1981).
\bibitem{M3} Mukai, S.: Moduli of vector bundles on K3 surfaces and symplectic manifolds. (Japanese) Sugaku Expositions \textbf{1}, 2, 139--174 (1988). Sugaku \textbf{39}, 3, 216--235 (1987).
\bibitem{M4} Mukai, S.: On the moduli space of bundles on K3 surfaces. I. Vector bundles on algebraic varieties (Bombay, 1984), 341--413, Tata Inst. Fund. Res. Stud. Math. \textbf{11}, Tata Inst. Fund. Res., Bombay (1987).
\bibitem{OG1} O'Grady, K.: Desingularized Moduli Spaces of Sheaves on a K3 Surface. J. Reine. Angew. Math. \textbf{512}, 49--117 (1999).
\bibitem{OG2} O'Grady, K.: A new six-dimensional irreducible symplectic variety. J. Algebraic Geometry \textbf{12}, 3, 435--505 (2003).
\bibitem{Pe} Perego, A.: The 2-factoriality of the O'Grady moduli spaces. Math. Ann. \textbf{346} (2), 367--391 (2009).
\bibitem{PR} Perego, A., Rapagnetta, A.: Deformations of the O'Grady moduli spaces. arxiv:1008.0190.
\bibitem{R1} Rapagnetta, A.: On the Beauville Form of the Known Irreducible Symplectic Varieties. Math. Ann. \textbf{321}, 77--95 (2008).
\bibitem{R2} Rapagnetta, A.: Topological invariants of O'Grady's six dimensional irreducible symplectic variety. Math. Z. \textbf{256}, 1--34 (2007).
\bibitem{S} Sorger, C.: On Moduli of $G-$bundles over Curves for exceptional $G$. Ann. Sci. \'Ecole Norm. Sup., \textbf{32}, 4, 127--133 (1999).
\bibitem{Y1} Yoshioka, K.: The Picard group of the moduli space of stable sheaves on a ruled surface, J. Math. Kyoto Univ. \textbf{36}, 279--309 (1996).
\bibitem{Y2} Yoshioka, K.: Chamber structure of polarizations and the moduli of stable sheaves on a ruled surface, Internat. J. Math. \textbf{7}, 411--431 (1996).
\bibitem{Y3} Yoshioka, K.: Albanese map of moduli of stable sheaves on abelian surfaces, arXiv:math/9901013.
\bibitem{Y4} Yoshioka, K.: Some notes on the moduli of stable sheaves on elliptic surfaces, Nagoya Math. J. \textbf{154}, 73--102 (1999).
\bibitem{Y5} Yoshioka, K.: Moduli spaces of stable sheaves on abelian surfaces. arxiv:0009001v2.
\end{thebibliography}
\end{document}